\documentclass[12pt]{amsart}
\usepackage{amssymb,amsmath,amsfonts,latexsym}
\usepackage{mathtools}
\usepackage{tikz-cd}
\usepackage[colorlinks=true, linkcolor=blue, citecolor=red, anchorcolor=orange,]{hyperref}

\newcommand{\newsection}[1]{\setcounter{equation}{0} \section{#1}}
\def\textmatrix#1&#2\\#3&#4\\{\bigl({#1 \atop #3}\ {#2 \atop #4}\bigr)}
\def\dispmatrix#1&#2\\#3&#4\\{\left({#1 \atop #3}\ {#2 \atop #4}\right)}
\newcommand{\be}{\begin{equation}}
	\newcommand{\ee}{\end{equation}}
\newcommand{\ben}{\begin{eqnarray*}}
	\newcommand{\een}{\end{eqnarray*}}

\newcommand{\bi}{\begin{itemize}}
	\newcommand{\ei}{\end{itemize}}

\newtheorem{theorem}{Theorem}[section]
\newtheorem{lemma}[theorem]{Lemma}
\newtheorem*{main theorem}{Main Theorem}{}
\newtheorem{proposition}[theorem]{Proposition}
\newtheorem{corollary}[theorem]{Corollary}
\theoremstyle{definition}
\newtheorem{definition}[theorem]{Definition}
\newtheorem{example}[theorem]{Example}

\newtheorem{conjecture}[theorem]{Conjecture}

\newtheorem{question}[theorem]{Question}
\newtheorem{remark}[theorem]{Remark}
\numberwithin{equation}{section}

\begin{document}

\title{On Solvability of Automorphism groups of Commutative algebras}

\author[Das]{Dibyendu Das}
\address{Indian Statistical Institute, Statistics and Mathematics Unit, 8th Mile, Mysore Road, Bangalore, 560059,
India}
\email{ddas4282@gmail.com}


	
\begin{abstract}
Let $A$ be a finite-dimensional commutative associative algebra with unity over an algebraically closed field $\mathbb{K}$. The purpose of the paper is to study the solvability of $G_A$, where $G_A$ is the identity component of $\text{Aut}_\mathbb{K}(A)$. Inspired by the work of Pollack (\cite{Pollack}), Saor\'in and Asensio have started this study for a commutative associative algebra $A$ when $\text{dim}(R/R^2)=2$, where $R$ is the Jacobson radical of $A$ (\cite{AS1}). In this paper, we give new sufficient conditions on $A$ so that $G_A$ is solvable without any restriction on $\text{dim}(R/R^2)$.
\end{abstract}
	
\maketitle

\tableofcontents

\newsection{Introduction and Main result}\label{sec: intro}

\noindent

Let $(V, 0)$ be an isolated singularity in $(\mathbb{C}^{n}, 0)$ defined by the zero set of a holomorphic function $f$. The moduli algebra $A(V)$ of $(V, 0)$ is $\mathbb{C}\{X_1,X_2,\dots,X_n\}/\langle f,\frac{\partial f}{\partial X_1},\dots,\frac{\partial f}{\partial X_n}\rangle$, where $\mathbb{C}\{X_1, X_2,\dots,X_n\}$ denotes the convergent power series ring. It is easy to see that $A(V)$ is a finite-dimensional $\mathbb{C}$-algebra and is an invariant of $(V, 0)$. The algebra of polynomials in $n$ indeterminates over $\mathbb{C}$ is denoted by $\mathbb{C}[X_1,\dots, X_n]$.  An \emph{isolated hypersurface singularity} is denoted by $V(f)=\{x\in \mathbb{C}^n:f(x)=0\}$ (in short IHS). By finite determinacy theorem (\cite{Greuel}, Section $2.2$ ), we may assume that a polynomial $f\in \mathbb{C}[X_1,\dots, X_n]$ defines an IHS at the origin, that is the hypersurface $V(f)$ has an isolated singularity at the origin. Let $J(f)=\langle\frac{\partial f}{\partial X_1},\dots,\frac{\partial f}{\partial X_n}\rangle$ and $A(V)$ be the quotient algebra $\mathbb{C}\{X_1,X_2\dots,X_n\}/\langle f,J(f)\rangle$. It is well known from (\cite{Greuel}, Lemma $2.3$) that $A(V)$ is a finite-dimensional algebra if and only if the hypersurface $V(f)$ has an isolated singularity at the origin. This algebra is called the \emph{moduli algebra} of the IHS. Mather and Yau have proved in \cite{Mather} that two IHS are biholomorphic equivalent if and only if their moduli algebras are isomorphic. Thus, the finite-dimensional local algebra $A(V)$ defines the IHS up to an analytic isomorphism. The \emph{Yau algebra} $L(V)\coloneq \text{Der}(A(V))$ is an important Lie algebra for understanding the IHS. In \cite{YAU2} and \cite{Yau}, it has been established by Yau that, if \( L(V) \) corresponds to the isolated hypersurface singularity defined by \( V(f) \), then \( L(V) \) is a solvable Lie algebra. By the above discussion, we have the following map from the space of isolated hypersurface singularities to the space of finite-dimensional associative commutative local algebras.

\vspace{0.1in}
\[
\begin{array}{ccc}
\left\{\begin{array}{c}
\text{Isolated hypersurface singu-} \\
\text{larities of dimension } n\\
 \text{ upto biholomorphic equivalence }
\end{array}\right\} & \longrightarrow & \left\{\begin{array}{c}
\text{Finite-dimensional associative } \\
\text{commutative local algebras}\\
\text{upto isomorphism }
\end{array}\right\} \\
\rotatebox{90}{$\in$} && \rotatebox{90}{$\in$} \\
(V,0) & \longrightarrow & A(V), \text{ the moduli algebra of } V
\end{array}
\]

\vskip2mm
\noindent
From Mather and Yau's theorem, it is clear that the above map is injective. However, the image of the above map is not understood. So, one may ask the following question: 

\vskip2mm
\begin{question}
Which finite-dimensional local commutative associative algebras can occur as the moduli algebra of some isolated hypersurface singularity?
\end{question}

\vskip2mm
The recognition problem constitutes a significant challenge within this field of study. Solvability of the derivation algebra of a moduli algebra by Yau gives a necessary condition for an algebra to be in the image of the above map. We already know that $\text{Aut}_{\mathbb{K}}(A)^{0}$ is solvable if and only if $\text{Der}_{\mathbb{K}}(A)$ is solvable over an algebraically closed field $\mathbb{K}$ with characteristic $0$. Therefore, studying the solvability of $\text{Aut}_{\mathbb{K}}(A)$ is important, where $A$ is a finite-dimensional associative commutative local algebra over $\mathbb{K}$. In this paper, we study the solvability of $G_A$ (the identity component of $\text{Aut}_{\mathbb{K}}(A)$) for an associative commutative algebra $A$ over an algebraically closed field $\mathbb{K}$. Our results are given for the case when $A$ is local over $\mathbb{K}$, to which everything can essentially be reduced (see Section \ref{sec: Solvability}). 

\vspace{0.1in}
The main inspiration for our work is Pollack's paper \cite{Pollack}, where he has proved different properties of $G_A$ from the properties of $A$. He gave a necessary and sufficient condition for the nilpotency of $G_A$, where $A$ is a finite-dimensional local associative algebra over an algebraically closed field $\mathbb{K}$. Subsequently, Saor\'in et al. established the conditions under which the solvability problem for local commutative algebras can be addressed, specifically when the dimension of \({R/R^2} \) equals $2$, where $R$ is the Jacobson radical of $A$ (refer to \cite{AS1}). In \cite{AS2}, these results have been refined for monomial algebras. Halperin's conjecture (see Section \ref{sec: Preliminaries}, Conjecture \ref{Conjecture 2.6}) is also an important problem in this context. Kraft and Procesi have solved a weaker version of the conjecture in \cite{kraft}. Perepechko has significantly contributed to understanding the solvability of $G_A$ (\cite{Perepechko}). He has solved Halperin's conjecture \ref{Conjecture 2.6} in \cite{Perepechko}. The problem of the solvability of the first Hochschild cohomology of a finite-dimensional associative algebra has been studied by many mathematicians, and it has a direct connection to this problem. For more details, we refer to the articles \cite{FR}, \cite{CSS}, \cite{DSSS}, \cite{S}, and the references therein. In Section $1$ of \cite{Pollack}, Pollack has shown that to study the solvability of $G_A$ for any finite-dimensional unital associative algebra $A$ over an algebraically closed field $\mathbb{K}$, it is enough to study the solvability of $G_A$ for a local associative algebra $A$ over $\mathbb{K}$. It has been proved in Corollary \ref{Corollary 3.14} that for a local finite-dimensional associative algebra $A$ over $\mathbb{K}$, the group $G_A$ is solvable, if the group $G_{G_r(A)}$ is solvable, where $G_{r}(A)$ denotes the associated graded by the radical algebra of $A$ (see Definition \ref{Definition 3.10}) and $G_{G_r(A)}$ is the identity component of $\text{Aut}_{\mathbb{K}}(G_{r}(A))$. Therefore, it is sufficient to investigate the solvability of $G_A$ for a local graded by the radical algebra $A$. In this article, we present new sufficient conditions for a local commutative associative algebra $A$ which is graded by the radical, ensuring that $G_A$ is solvable (see Theorem \ref{Theorem 5.1}). More precisely, we have the following main result without any restriction on $\text{dim}(R/R^2)$:

\vspace{0.2in}
\begin{main theorem}
    Let $A$ be a local commutative graded by the radical algebra over an algebraically closed field $\mathbb{K}$ with $\emph{dim}(R/R^2)=n$. Then by quiver and relation representation of $A$, we have an isomorphism $A\cong \mathbb{K}[X_1,\dots, X_n]/I$ with $\langle X_1,\dots, X_n\rangle^l\subsetneqq I\subsetneqq \langle X_1,\dots, X_n\rangle^2$, where $l>2$ is the Lowey length of $A$, and $I$ is a homogeneous ideal. Let $W$ be the minimal degree subspace associated to $A$, i.e., the collection of all non-zero homogeneous polynomials of minimal degree in $I$ and $0$ (see Remark \ref{Remark 4.18} ). Let any of the conditions listed below be satisfied:
    \vskip2mm
    \begin{enumerate}
    \vskip2mm
    \item $\emph{\text{char }}\mathbb{K}=0$ and $\emph{\text{dim }}W=n$ with $W$ having an ordered basis, whose elements form a regular sequence in the polynomial ring $\mathbb{K}[X_1,\dots, X_n]$;
    \vskip2mm
    \item $\emph{\text{char }}\mathbb{K}=0$ and $W$ contains a non-singular homogeneous polynomial of degree at least $3$;
\item $\emph{\text{char }}\mathbb{K}=p>3$ and $W$ contains a non-singular homogeneous polynomial of degree $d$ with $2<d<p$.  
    \end{enumerate}
    \vskip2mm
 Then, $G_A$ is solvable.
\end{main theorem}

\vskip2mm
\noindent
 Let $\mathbb{K}$ be an algebraically closed field and $A$ be a finite-dimensional associative unital algebra over the field $\mathbb{K}$. Let us denote $\text{Aut}_{\mathbb{K}}(A)$, the group of all $\mathbb{K}$-algebra automorphisms of $A$. Then $\text{Aut}_{\mathbb{K}}(A)$ is an affine algebraic group and $G_A$ denotes its identity component. Associative algebras and affine algebraic groups have a rich structure theory. From now onwards, algebra means a finite-dimensional associative unital algebra. Both (\cite{Humphreys}) and (\cite{Borel}) serve as general references for algebraic groups, and we use the theory of associative algebras from (\cite{Pierce}, \cite{Auslander}, \cite{Curtis}). The concept of representation of a split basic finite-dimensional associative algebra by \emph{quiver and relation} is taken from (\cite{Auslander}, \cite{PG}, \cite{ISD}). A group $G$ over a perfect field $K$ is solvable if and only if $G\otimes \overline{K}$ is solvable over $\overline{K}$, where $\overline{K}$ is an algebraic closure of $K$. If $K$ is an algebraically closed field, then for any local finite-dimensional associative algebra $A$ over $K$, we have $A/R\cong K$ (split local associative algebra), where $R$ is the Jacobson radical of $A$. The results we prove here are generally true for any split local finite-dimensional associative commutative algebra over a perfect field $K$, i.e., if $A/R\cong K$. For a split local associative algebra $A$ over a perfect field $K$, $G_A$ is solvable over $K$ if and only if $G_{A\otimes \overline{K}}$ is solvable over $\overline{K}$. For simplicity, we prove all our results over an algebraically closed field. This article will deal with finite-dimensional associative algebras over an algebraically closed field.
\vskip2mm
\noindent
 The rest of the paper is organized as follows: In Section \ref{sec: Preliminaries}, we recall some well-known results on the structure of finite-dimensional associative algebras and provide some basic definitions that are used throughout the paper. Section \ref{sec: Solvability} reviews the result of Pollack and some well-known facts about the solvability of the automorphism groups of different classes of associative algebras. Section \ref{sec: commutative algebra} is devoted to the study of the automorphism groups of finite-dimensional local commutative algebras by using quiver and relation representation. In Section \ref{sec: Main}, we establish the main result (Theorem \ref{Theorem 5.1}) of this article. Subsequently, we illustrate a few examples of new classes of algebras that exhibit $G_A$ is solvable. Finally, in Section \ref{sec: remark}, we provide examples that demonstrate the significance of the various conditions in the main theorem.
 
\section{Preliminaries}\label{sec: Preliminaries}
\vskip2mm
\noindent
This section discusses the structure theory of finite-dimensional associative algebras and provides some basic definitions that are necessary for subsequent sections.
\vskip2mm
\noindent
Throughout this paper, $A$ is a finite-dimensional associative algebra over a field $\mathbb{K}$ with unity, and $R$ is its Jacobson radical (in short, we call it radical), where $\mathbb{K}$ is an algebraically closed field. Unless specified otherwise, there are no other restrictions on the field $\mathbb{K}$. For any affine algebraic group $G$ over $\mathbb{K}$, $G^{0}$ denotes its identity component.

\vskip2mm
\noindent
Next, we recall the principal structure theorem for associative algebras: the Wedderburn-Malcev decomposition.

\vskip2mm
\begin{theorem}\label{Theorem 2.1}
\textbf{$($\emph{Wedderburn-Malcev}}, \cite{Curtis}, \emph{Theorem} $72.19)$ Let $A$ be a finite-dimensional associative algebra over $K$ (arbitrary field) with radical $R$ such that the residue class algebra $B=A/R$ over $K$ is separable. Then, there exists a semisimple subalgebra $A_s$ of $A$ such that $ A=A_s\oplus R$ (vector space decomposition). If there exists a semisimple subalgebra $A'_s$ such that $A=A'_s\oplus R$, then there exists an element $r\in R$ such that $A'_s=(1+r)A_s(1+r)^{-1}$.
\end{theorem}

\vskip2mm

\begin{corollary}\label{Corollary 2.2}
$($\cite{Pierce}, \emph{Page} $211)$ If $K$ is a perfect field, and $A$ is finite-dimensional associative $K$-algebra, then there is a semisimple subalgebra $A_s$ of $A$ such that $A/R\cong A_s$ and $A=A_s\oplus R$. Moreover, $A_s$ is unique upto
conjugation by units of the form $(1+r)$, where $r\in R$.
\end{corollary}

\vspace{0.1in}
\noindent
Thus, in the context of an algebraically closed field \( \mathbb{K} \), every finite-dimensional associative algebra \( A \) can be decomposed as \( A = A_s \oplus R \) (vector space decomposition), as established by Theorem \ref{Theorem 2.1}. In this decomposition, \( A_s\cong A/R \) is a semisimple subalgebra that is isomorphic to \( \bigoplus_{i=1}^{m} M_{n_{i}}(\mathbb{K}) \), and \( R \) denotes the Jacobson radical of the algebra \( A \).
\vskip4mm
\begin{definition}\label{Definition 2.3}
A finite-dimensional associative algebra $A$ over a perfect field $K$ is called \emph{split algebra}, if $A/R\cong \bigoplus_{i=1}^{m}M_{n_i}(K)$, $m\in\mathbb{N}$.
\end{definition}
\vskip2mm
\begin{definition}\label{Definition 2.4}
Let $K$ be a perfect field. An algebra $A$ over $K$ is a \emph{split basic algebra}, if $A/R\cong \bigoplus_{i=1}^{n} K$, where $n\geq 1$. In particular, if $A$ is an algebra over $K$ such that $A/R\cong K$, then we call $A$ a \emph{split local algebra} over $K$. Over an algebraically closed field $\mathbb{K}$, an associative algebra $A$ is a basic (resp. local) algebra, if $A/R\cong \bigoplus_{i=1}^{n}\mathbb{K}$ (resp. $A/R\cong \mathbb{K})$.  
\end{definition}
\vskip2mm
\begin{definition}\label{Definition 2.5}
Let $K$ be a perfect field. An ordered sequence of polynomials $\{f_1,\dots,f_m\}$ in $K[X_1,\dots,X_n]$ with $1\leq m\leq n$ is said to be a regular sequence in $K[X_1,\dots,X_n]$ if both of the following conditions hold:
\begin{enumerate}
\item The ideal $\langle f_1,\dots,f_m\rangle \neq K[X_1,\dots,X_n].$
\vskip2mm
\item For all $i$, $0\leq i\leq m-1$, $f_{i+1}$ is not a zero divisor in $K[X_1,\dots,X_n]/\langle f_1,\dots,f_i\rangle$.
\end{enumerate}
\end{definition}
\vskip2mm
\begin{conjecture}$($Halperin, \cite{kraft}$)$\label{Conjecture 2.6}
Let $f_1,\dots,f_n$ be a regular sequence in the polynomial ring $\mathbb{C}[X_1,\dots, X_n]$. Then the connected component of the automorphism group of the finite-dimensional algebra \[\mathbb{C}[X_1\dots, X_n]/\langle f_1,\dots,f_n\rangle\] is solvable.

\end{conjecture}
\vskip2mm
\begin{definition}\label{Definition 2.7}
$($\cite{SK}$)$
Let $G$ be a connected linear algebraic group, and $\rho$ a rational representation of $G$ on a finite-dimensional vector space $V$, all defined over the perfect field $K$. We call such a triple $(G, \rho, V)$ a prehomogeneous vector space if $V$ has a Zariski-dense open $G$-orbit.
\end{definition}
\vskip 2mm
\noindent
To understand the automorphism group of an associative finite-dimensional algebra $A$, we use the quiver and relation representation of $A$. The representation of an associative algebra $A$ over an algebraically closed field $\mathbb{K}$ using a quiver has a long tradition. Let $\Gamma(A)$ be the quiver of $A$ and $V(\Gamma)$ (resp. $A(\Gamma))$ be the set of vertices (resp. arrows). It is well-known that $A$ is Morita equivalent to the quotient algebra $\mathbb{K}[\Gamma]/I$, where $\mathbb{K}[\Gamma]$ is the path algebra and $J^m \subset I\subset J^2$, where $J$ denotes the ideal of $\mathbb{K}[\Gamma]$ consisting of all linear combinations of paths of length $\geq1$, and $I$ is called the \emph{adequate ideal} or \emph{admissible ideal} for $A$. If $A$ is a basic algebra over an algebraically closed field $\mathbb{K}$, then $A\cong \mathbb{K}[\Gamma]/I$ (see \cite{PG}, \cite{ISD}, \cite{Auslander}).

\section{Solvability of $G_A$}\label{sec: Solvability}
\vskip2mm
\noindent
In this section, we review Pollack's result related to the solvability of $G_A$ and recall the well-known facts about the solvability of the automorphism groups of finite-dimensional associative algebras.
\vskip2mm
\noindent
Let $A^{\times}$ denote the group of units in $A$ and $\mathbb{K}^{\times}=\mathbb{K}-\{0\}$. We have an obvious map from $A^{\times}$ to $\text{Aut}_{\mathbb{K}}(A)$ given  by $\psi:A^{\times}\rightarrow \text{Aut}_{\mathbb{K}}(A) \text{ with }\psi(a)=\text{Int}(a)$, where $\text{Int}(a)$ is the automorphism  $x\mapsto axa^{-1}$.

\vskip2mm
\noindent
This is a morphism of algebraic groups, whose image is $\text{Inn}(A)$, the group of inner automorphisms of $A$. Observe that $\{1+r:r\in R\}$ is a closed, connected, unipotent subgroup of $A^{\times}$, whose image $\hat{R}\coloneq \{\text{Int}(1+r):r\in R\}$ is normal in $\text{Aut}_{\mathbb{K}}(A)$, and $\hat{R}\subset U$, where $U$ is the unipotent radical of $\text{Aut}_{\mathbb{K}}(A)$.
\vskip2mm
\noindent
Now we will decompose the identity component of the automorphism group of finite-dimensional associative algebras, similar to Pollack's paper \cite{Pollack}. The following result is clear from Pollack's paper; we include the proof here for the sake of completeness.
\vskip2mm
\begin{theorem}$($\cite{Pollack}$)$\label{Theorem 3.1}
    Let $A$ be a finite-dimensional associative algebra over an algebraically closed field $\mathbb{K}$, then $G_A$ is solvable if and only if $G_{A, A_s}\coloneq \{\sigma\in G_A:\sigma|_{A_s}=Id\}$ is solvable and $A$ is a basic algebra, where $A_s$ is a semisimple subalgebra $A$ such that $A=A_s\oplus R$ (vector space decomposition).
\end{theorem}
\vskip1mm
\begin{proof}
First, we will prove some lemmas to decompose the group $G_A$.  Let $A= A_s\oplus R$ be a Wedderburn-Malcev decomposition of $A$, where $A_s\cong A/R$ is a semisimple subalgebra of $A$.
\begin{lemma}\label{Lemma 3.2}
$G_A=\hat{R}G_{A}^{A_s}$, where $\hat{R}=\{\emph{\text{Int}}(1+r)\in G_A : r\in R\}$, $R$ is the radical of $A$ and $G_{A}^{A_s}\coloneq\{\sigma \in G_A: \sigma(A_s)=A_s\}$.
\end{lemma}
\vskip1mm
\begin{proof}
Let $\sigma \in G_A$, then $A_s\oplus R=A=\sigma(A)=\sigma(A_s)\oplus R$. Therefore, using Corollary \ref{Corollary 2.2} we have $\sigma(A_s)=\rho(A_s)$, for some $\rho \in \hat{R}$. Thus, $\rho^{-1}\sigma\in G_{A}^{A_{s}}$ and it follows that $G_A=\hat{R}G_{A}^{A_s}$. 
\end{proof}

\vskip2mm
\begin{lemma}\label{Lemma 3.3}
$G_{A}^{A_s}\cong \emph{Inn}(A_{s}^{\times})G_{A,A_s}$, where $\emph{Inn}(A_s^{\times})\coloneq\{\emph{Int}(x)\in G_{A/R} : x\in A_{s}^{\times}\}$, and $G_{A,A_s}\coloneq\{\sigma\in G_A: \sigma|_{A_s}=Id\}$.
\end{lemma}
\vskip1mm
\begin{proof}
We have the induced map 
\begin{equation*}\label{eq11}
G_A\rightarrow \text{Aut}_{\mathbb{K}}(A/R)
\end{equation*}
whose image is closed and connected so contained in $\text{Inn}(A/R)=G_{A/R}$ as $A/R \cong \bigoplus_{i}M_{n_i}(\mathbb{K})$. Therefore, we have an exact sequence 
\begin{center}
\begin{tikzcd}
1\arrow[r] & G_{A,A_s} \arrow[r] &G_{A}^{A_s} \arrow[r,"res"] & \text{Inn}(A_s^{\times})\arrow[bend left=33]{l}{s}\arrow[r] & 1,
\end{tikzcd}
\end{center}
\vskip1mm
\noindent
where $s$ is a section $s:\text{Inn}(A_s^{\times})\rightarrow G_{A}^{A_s}$ defined by $s(\text{Int}(a)_{|A_s})=\text{Int}(a)$, $a\in A_s^{\times}$, $res(\phi)=\phi_{|A_s}$, $\phi\in G_{A}^{A_s}$ with $(res)\circ s=Id$. Therefore, the sequence is split exact. Hence, the result follows.
\end{proof}
\vskip2mm
\noindent
Consider the map $\text{Int}:A_s^{\times}\rightarrow \text{Inn}(A_s^{\times})$ given by $\text{Int}(a)(x)=axa^{-1}$, which is surjective and $\text{Ker}(\text{Int})\subseteq Z(A_s^{\times})$, $Z(A_s^{\times})$ denotes the center of $A_s^{\times}$.
\vskip4mm
\noindent
Using the previous two lemmas, we conclude that $G_A$ is solvable if and only if $\text{Inn}(A_s^{\times})$ and $G_{A,A_s}$ are solvable. However, $A_s\cong \bigoplus_{i} M_{n_{i}}(\mathbb{K})$, hence $A_s^{\times}\cong \bigoplus_i \text{GL}_{n_i}(\mathbb{K})$ and $Z(A_s^{\times})\cong \bigoplus_i \mathbb{G}_{m,n_i}$, where $\mathbb{G}_{m,n_i}$ denotes the torus of scalar matrices in $\text{GL}_{n_i}(\mathbb{K})$. Therefore, $\text{Ker}(\text{Int})$ is solvable. Hence, $G_A$ is solvable if and only if $A_s^{\times}$ and $G_{A,A_s}$ are solvable. Therefore, $G_A$ is solvable if and only if $A_s$ is a direct sum of fields and $G_{A,A_s}$ is solvable. The result now follows.
 \end{proof}
\vskip2mm
\begin{remark}\label{Remark 3.4}
The above theorem is true for any split finite-dimensional algebra $A$ over a perfect field $K$, i.e., if $A/R\cong \bigoplus_iM_{n_i}(K)$.
\end{remark}
 \vspace{0.1in}
 \noindent
 The following example shows that $A$ is a basic algebra, which is a necessary condition for the solvability of $G_A$, but not sufficient.
 \vskip2mm
\begin{example}\label{Example 3.5}
 Let $A=\mathbb{K}[X_1,X_2,\dots,X_n]/\langle X_1,X_2,\dots,X_n\rangle^l$ for some $l\geq 2$ and $n\geq 2$, then $\text{Aut}_\mathbb{K}(A)\cong \text{GL}_n(\mathbb{K})\ltimes U$, where $U$ is a unipotent group (see \cite{AS1}, Theorem $1.3$). Therefore, $G_A$ is not solvable, though $A/R\cong \mathbb{K}$.  Note that $A_s\cong \mathbb{K}$ and $G_A=G_{A,A_s}=\text{GL}_n(\mathbb{K})\ltimes U$, which is not solvable.
\end{example}

\vspace{0.1in}
\noindent 
Since for the solvability of $G_A$, $A$ is necessarily a basic algebra, it follows that $A=A_s\oplus R\cong \bigoplus_i A_i$, where $A_i$'s are local algebras as $A_s\cong \bigoplus_i \mathbb{K}$. It implies $G_A\cong \prod _i G_{A_i}$. Therefore, it is enough to study the automorphism group of finite-dimensional local associative algebras to study the solvability of $G_A$ in general.

\vspace{0.1in}
\noindent
The following theorems are the known for the solvability of $G_A$ over an algebraically closed field $\mathbb{K}$:
\vskip2mm
\begin{theorem}\label{Theorem 3.6}
$($\cite{Pollack}, \emph{Theorem} $1.6)$
Let $A$ be a local associative algebra over the field $\mathbb{K}$ and assume that $G_A\neq \{Id\}$ is nilpotent. Then $G_A$ is either a torus or a unipotent group.
\vskip2mm
\noindent
$(a)$The following are equivalent for a local associative algebra $A$ with $G_A$ nilpotent;
\begin{enumerate}
\item $G_A$ is a non-trivial torus.
\vskip2mm
\item $R^2=0$ and $R\neq 0$.
\vskip2mm
\item $A\cong \mathbb{K}[X]/\langle X\rangle^2$
\end{enumerate}
\vskip2mm
\noindent
$(b)$The following are equivalent for a local associative algebra $A$ with $G_A$ is nilpotent;
\begin{enumerate}
\item $G_A$ is unipotent.
\vskip2mm
\item $R^2\neq 0$.
\vskip2mm
\item $\emph{dim}(R)>1$.
\end{enumerate}
\end{theorem}
\vspace{0.1in}
\begin{theorem}\label{Theorem 3.7}
$($\cite{Perepechko}, \emph{Corollary} $4.3)$
Let $\mathbb{K}$ be an algebraically closed field with characteristic $0$. Suppose the ideal $I\subset \mathbb{K}[X_1,\dots, X_n]$ with $m$ generators and an integer $l>1$ is such that each of the following conditions holds.
\begin{enumerate}
\vskip2mm
\item The quotient algebra $A=\mathbb{K}[X_1,\dots X_n]/I$ is finite-dimensional.
\vskip2mm
\item For any maximal ideal $\mathfrak{M}\subset \mathbb{K}[X_1,\dots,X_n]$, either $I\not\subseteq \mathfrak{M}$ or $I\subseteq \mathfrak{M}^{l}$.
\vskip2mm
 \item $m<n+l-1$.
\end{enumerate}
\vskip2mm
Then, the identity component $\emph{\text{Aut}}_{\mathbb{K}}(A)^{0}$
is solvable.
\end{theorem}
\vspace{0.1in}
\begin{theorem}\label{Theorem 3.8}
$($\cite{AS1}, \emph{Theorem} $2.6)$ Let $A$ be a local commutative associative algebra over the algebraically closed field $\mathbb{K}$ with $\emph{\text{char }}\mathbb{K}=0$. If $\emph{dim}(R/R^2)=2$, then either $A\cong \mathbb{K}[X,Y]/\langle X,Y\rangle^l$, where $l\geq 2$, or $G_A$ is solvable.    
\end{theorem}
\vspace{0.1in}
\noindent
Now we define some terminology to state the next theorem.
\vskip2mm
\noindent
Let $V(\Gamma)$ be the set of vertices for the quiver $\Gamma$ and $\overline{V(\Gamma)}^{2}$ denotes the tuples $(v,w)$ such that $A(v,w)\neq \emptyset$, where $A(v,w)$ is the set of arrows from $v$ to $w$. Let $\mathbb{K}[\Gamma]$ be the path algebra and $I$ be its admissible ideal. As per the notation in \cite{AS2}, $\alpha \cong_I\beta$ means for every path $p\in I$ in which $\beta$ occurs, the replacement once of $\beta$ by $\alpha$ still yields a path in $I$, similarly for path $p'\in I$ in which $\alpha$ occurs, the replacement once of $\alpha$ by $\beta$ again yields a path in $I$.

\vskip2mm
\begin{theorem}\label{Theorem 3.9}
$($\cite{AS2}, \emph{Corollary} $2.22)$
Let $\Gamma$ be a finite quiver and $I$ a monomial adequate ideal of $\mathbb{K}[\Gamma]$, where $\mathbb{K}$ is an algebraically closed field with $\emph{\text{char }}\mathbb{K}=0$. The following assertions are equivalent for the algebra $A\cong \mathbb{K}[\Gamma]/I$: 
\begin{enumerate}
\vskip2mm
\item $G_A$ is a solvable group.
\vskip2mm
\item For each $(v,w)\in \overline{V(\Gamma)}^2$, if $\alpha,\beta\in A(v,w)$ and $\alpha\cong_I\beta$ then $\alpha=\beta$.
\end{enumerate}
\end{theorem}

\vskip2mm

\begin{definition}\label{Definition 3.10}
A grading on an algebra $A$ is a family of vector subspaces $(A_n)_{n\in\mathbb{Z}}$ of $A$ such that $A=\bigoplus_{n\in\mathbb{Z}} A_n$ and $A_n.A_m\subseteq A_{n+m}$ for all $n,m\in \mathbb{Z}$. The grading is called \emph{by the radical}, when $A_m=0$ for $m<0$ and $R^m=\bigoplus_{n
\geq m} A_n$, for each $m\geq 0$ (see \cite{Saorin}, Proposition $1.1$ for equivalent definitions$)$.  For any associative algebra $A$ with $R^{0}=A$, we can associate a graded by the radical algebra, $G_{r}(A)\coloneq \bigoplus_{m\geq 0} R^m/R^{m+1}$ (\emph{associated graded by the radical algebra}), as a vector space, and where the multiplication $(a+R^{m+1}). (b+R^{n+1})=(ab+R^{m+n+1})$ for $a\in R^{m}$ and $b\in R^{n}$, extends to all of $G_{r}(A)$ by linearity. An algebra $A$ is said to be \emph{graded by the radical algebra} if it has a grading by the radical, or equivalently, $A$ is isomorphic to $G_{r}(A)$. For more details, we refer to \cite{AS1} and \cite{Saorin}.
\end{definition}
\vskip2mm
\noindent
\begin{remark}\label{Remark 3.11}
It is easy to see that there are well-defined maps for a local algebra $A$ over an algebraically closed field $\mathbb{K}$, which are given below. Since for any automorphism $\phi\in \mathrm{Aut}_{\mathbb{K}}(A)$, it induces an automorphism $\bar{\phi}:R/R^2\rightarrow R/R^2$. We have a morphism of algebraic groups $\Phi_A: \mathrm{Aut}_{\mathbb{K}}(A)\rightarrow \mathrm{GL}(R/R^2)$, where $\Phi_A(\phi)=\bar{\phi}$. If $\phi\in \mathrm{Aut}_{\mathbb{K}}(A)$, then we can define a $\mathbb{K}$-linear bijective map $\Psi_A(\phi)$, which is just the direct sum of the maps induced by $\phi$ on all $R^m/R^{m+1}$, $m\in \mathbb{N}\cup\{0\}$. It is clear that $\Psi_A(\phi)$ is actually an algebra automorphism, thus yields a map $\Psi_A:\mathrm{Aut}_{\mathbb{K}}(A)\rightarrow \mathrm{Aut}_{\mathbb{K}}(G_{r}(A))$. Similarly, for every $m\in\{2,\dots,n\}$, where $R^{n+1}=0$, we have a canonical map $\phi_{\mathbb{K}}: \mathrm{Aut}_{\mathbb{K}}(A)\rightarrow \mathrm{Aut}_{\mathbb{K}}(A/R^m)$ that takes $\phi\in \mathrm{Aut}_{\mathbb{K}}(A)$ to the induced automorphism $\tilde{\phi}: A/R^m\rightarrow A/R^m$.  It is easy to see that we have similar kind of maps $\Phi_{G_{r}(A)}: \mathrm{Aut}_{\mathbb{K}}(G_{r}(A))\rightarrow \mathrm{GL}(R/R^2)$ and $\Phi_{A/R^{m}}:\mathrm{Aut}_{\mathbb{K}}(A/R^m)\rightarrow \mathrm{GL}(R/R^2)$.  All of the above maps $(\Phi_A, \Psi_A, \phi_{\mathbb{K}},\Phi_{G_{r}(A)}, \Phi_{A/R^m}$) are morphism of algebraic groups.
\end{remark}
\vskip2mm
\noindent
 The following result is known from \cite{Pollack}. We include the proof here for convenience.
\vskip2mm
\begin{proposition}\label{Proposition 3.12}
Let $A$ be a local associative algebra over $\mathbb{K}$ and $G_r(A)$ is defined as above, i.e., $G_{r}(A)\coloneq (\mathbb{K}\oplus\bigoplus_{m=1}^{n}(R^m/R^{m+1}))$ be the associated graded by the radical algebra of $A$. Consider the following commutative diagrams $(1)$ and $(2)$:
\begin{equation}\label{eq6}
\begin{minipage}{0.45\linewidth}
\begin{tikzcd}
    \emph{\text{Aut}}_{\mathbb{K}}(A) \arrow[dr, "\Phi_A"] \arrow[d, "\Psi_A"']\\
    \emph{\text{Aut}}_{\mathbb{K}}(G_{r}(A)) \arrow[r,"\Phi_{G_{r}(A)}"'] & \emph{GL}(R/R^2)
\end{tikzcd}
\begin{center}
\textbf{(1)}
\end{center}
\end{minipage}
\begin{minipage}{0.45\linewidth}
\begin{tikzcd}
   \emph{ \text{Aut}}_{\mathbb{K}}(A) \arrow[dr, "\Phi_A"] \arrow[d, "\phi_{\mathbb{K}}"']\\
   \emph{ \text{Aut}}_{\mathbb{K}}(A/R^m) \arrow[r,"\Phi_{A/R^m}"'] & \emph{GL}(R/R^2)
\end{tikzcd}
\begin{center}
\textbf{(2)}
\end{center}
\end{minipage}
\end{equation}
Then each of the induced morphisms $(\Phi_A, \Phi_{G_{r}(A)}, \Phi_{A/R^m})$ has a unipotent kernel. Therefore, the rank of $\emph{\text{Aut}}_{\mathbb{K}}(A)$ is bounded by the dimension of $R/R^2$.
\end{proposition}
\vskip1mm
\begin{proof}
We will prove that $\text{Ker}(\Phi_A)$ is unipotent; other cases follow similarly. Since $A$ is a local algebra over $\mathbb{K}$, $\text{Aut}_{\mathbb{K}}(A)=\text{Aut}_{\mathbb{K}}(R)$. Therefore, \[\text{Ker}(\Phi_A)=\{\phi\in \text{Aut}_{\mathbb{K}}(R):\phi(x)-x\in R^2, \forall x\in R\}=\{\phi\in \text{Aut}_{\mathbb{K}}(R):(\phi -Id)(R)\subset R^2\}.\] Let $x\in R$ and $(\phi-Id)(x)=\sum_{i=1}^n x_iy_i$, where each $x_i,y_i\in R$. Then,
\begin{align*}
(\phi-Id)\bigl((\phi -Id)(x)\bigl)={} & 
\sum_{i=1}^n\phi(x_iy_i)-\sum_{i=1}^nx_iy_i.
\end{align*}
It implies 
\begin{flalign*}
 (\phi-Id)^2(x) ={}&
  \!\begin{aligned}[t] 
    &\sum_{i=1}^{n}\phi(x_i)\phi(y_i)-x_1\phi(y_1)+x_1\phi(y_1)-\dots\\ 
     &-x_n\phi(y_n)+x_n\phi(y_n)-\sum_{i=1}^nx_iy_i
  \end{aligned}\\
  ={}&
  \!\begin{aligned}[t]
    &\sum_{i=1}^n(\phi-Id)(x_i)\phi(y_i)+\sum_{i=1}^nx_i(\phi-Id)(y_i) \\
  \end{aligned}\\
\end{flalign*}
\noindent
However, $(\phi-Id)(x)\in R^2, \forall x\in R$. Hence, $(\phi-Id)^2(R)\subset R^3.$ Suppose $R^m=0$ for some $m\geq 2$. So, we can iterate the process and show that $(\phi-Id)^{m-1}=0.$ Since the kernel is a unipotent group, any maximal torus of $\text{Aut}_{\mathbb{K}}(A)$ embeds in $\text{GL}(R/R^2)$. Hence, the rank of any maximal torus is less than or equal to $\text{dim}(R/R^2)$.
\end{proof}
 \vskip2mm
 \begin{remark}\label{Remark 3.13}
   Restricting $\Phi_A$ to $ G_A$ yields a homomorphism from $ G_A$ to $\text{GL} (R/R^2)$. By a similar line of argument, we can prove that $\text{Ker}(\Phi_A|_{G_{A}})$ is a unipotent group. Similarly, we can apply restrictions to other mappings to their identity components. Still, the proposition remains valid. To avoid notational ambiguity, the notation used for the identity component mappings will be the same as $\Phi_A$ and $\Psi_A$.
 \end{remark}
\vskip2mm
\begin{corollary}\label{Corollary 3.14}
  Let $A$ be a local associative algebra over $\mathbb{K}$. If $\emph{\text{Aut}}_{\mathbb{K}}(G_{r}(A))$ $($resp. $G_{G_{r}(A)})$ is solvable, then $\emph{\text{Aut}}_{\mathbb{K}}(A)$ $($resp. $G_A)$ is solvable, where $G_{G_{r}(A)}\coloneq \emph{\text{Aut}}_{\mathbb{K}}(G_{r}(A))^{0}$.
\end{corollary}
\vskip1mm
\begin{proof}
   Since $\text{Ker}(\Psi_A)$ is a unipotent group, hence solvable, therefore $G_A$ is solvable if and only if $\text{Im}(\Psi_A$) is solvable. However, $\text{Im}(\Psi_A) \subset \mathrm{Aut}_{\mathbb{K}}(G_{r}(A))$ (resp. $\mathrm{Im}(\Psi_A)\subset G_{G_{r}(A)}$). Hence, we are done. 
\end{proof}
\vskip2mm
\noindent
Therefore, we will mainly focus on graded by the radical local algebras. For more about grading by the radical associative local algebras, we refer to \cite{Saorin}.

\vskip4mm
\noindent
Let $A$ be an arbitrary basic not necessarily commutative algebra. If $\Gamma$ is the quiver of $A$ and $I$ is an ideal of $\mathbb{K}[\Gamma]$, $I_{*}$ will denote the ideal of $\mathbb{K}[\Gamma]$ generated by the homogeneous components with the smallest degree of elements of $I$. We state the following proposition for use in the future.

\vskip2mm

\begin{proposition}\label{Proposition 3.15}
$($\cite{AS1}, \emph{Proposition} $2.1)$ Let $A$ be a basic not necessarily commutative algebra, with quiver $\Gamma$ and $G_{r}(A)$ be the associated graded by the radical algebra of $A$. If $I$ is an adequate ideal for $A$ in $\mathbb{K}[\Gamma]$, then $I_{*}$ is an adequate ideal for $G_{r}(A)$. Conversely, if $L$ is a homogeneous adequate ideal for $G_{r}(A)$, then there is an adequate ideal $I$ for $A$ such that $I_{*}=L$.
\end{proposition}

\section{Finite Dimensional Commutative Algebras}\label{sec: commutative algebra}

\vskip2mm
\noindent
The purpose of this section is to understand the automorphism groups of finite-dimensional local commutative algebras by quiver and relations. From this point onward, our primary emphasis will be on the automorphism groups of finite-dimensional commutative associative local algebras with unity over $\mathbb{K}$, where $\mathbb{K}$ will always denote an algebraically closed field. Using the representation of finite-dimensional associative algebras by quiver and relation, one knows that any such algebra $A$ is isomorphic to $\mathbb{K}[X_1,\dots, X_n]/{I}$, where $I$ is an ideal of $\mathbb{K}[X_1,\dots,X_n]$ such that there exists an $l\geq2$ with $\langle X_1,\dots,X_n \rangle^l \subseteq I \subseteq \langle X_1,\dots,X_n \rangle^2$ (see \cite{AS1}, \cite{PG}). The minimum $l$ is called the Lowey length of $A$, denoted by $L(A)$. Every such ideal is called an \emph{admissible} or \emph{adequate ideal} for $A$ in $\mathbb{K}[X_1,\dots,X_n]$. A finite set of generators of an adequate ideal is called \emph{adequate sets of relations} for $A$. When $A$ is identified with $\mathbb{K}[X_1,\dots,X_n]/I$, $R^m=(\text{Jacobson}(A))^m$ gets identified with $[(X_1,\dots,X_n)^m+I]/I$ and $n = \text{dim}(R/R^2)$. Let $A$ be identified with $\mathbb{K}[X_1,\dots,X_n]/I$, where $I$ is an admissible ideal. Suppose we denote $x_i=X_i+I$, then $\{\bar{x_i}=x_i+R^2: i\in\{1,\dots,n\}\}$ is a basis of $R/R^2$. By taking matrices with respect to that basis, $\text{GL}(R/R^2)$ gets identified with $\text{GL}_n(\mathbb{K})$. The notion of change of variables is defined as follows.

\vskip2mm

\begin{definition}\label{Definition 4.1}
$($\cite{AS1}$)$ A homomorphism of algebras $F:\mathbb{K}[X_1,\dots,X_n]\rightarrow \mathbb{K}[X_1,\dots,X_n]$ is said to be a change of variables in $\mathbb{K}[X_1,\dots,X_n]$, if there is $n\times n$ matrix $(\lambda_{ij})\in\text{GL}_n(\mathbb{K})$ such that $F(X_j)=\sum_{1\leq i\leq n}\lambda_{ij}X_i$ modulo $(\langle X_1,\dots,X_n\rangle^2)$ for every $j=1,\dots,n$. In particular, if no monomial of degree $\geq 2$ appears in the change of variables $F$, i.e., $F(X_j)=\sum_{1\leq i\leq n}\lambda_{ij}X_i$ with $(\lambda_{ij})\in \text{GL}_n(\mathbb{K})$ for all $j\in\{1,\dots,n\}$, then we call such a change of variables a \emph{linear change of variables}, and $(\lambda_{ij})_{n\times n}$ is called the \emph{linear change of variables matrix} and denoted by $M_{F}=(\lambda_{ij})_{n\times n}$. For any linear change of variables $F$, we have $F(f)(X)=f(XM_{F})$, where $f\in \mathbb{K}[X_1,\dots,X_n]$ and $X=(X_1,\dots,X_n)$.
\end{definition}
\vskip2mm
\noindent
By using Theorem $1.5$ from \cite{Saorin} for a commutative local algebra $A$ over $\mathbb{K}$, we get that $A$ is graded by the radical if and only if $A$ admits a homogeneous admissible ideal $I$.

\vskip2mm

\begin{proposition}\label{Proposition 4.2}
$($\cite{AS1}$)$  Let $A$ be a local commutative algebra with Lowey length $l$ and $\emph{dim}(R/R^2) = n$. Suppose that $I, I'$ are two ideals of $\mathbb{K}[X_1,\dots,X_n]$, the first being adequate for $A$. The following assertions are equivalent for a homomorphism of algebras
\begin{equation*}\label{eq19}
\phi: \mathbb{K}[X_1,\dots,X_n]/I \longrightarrow \mathbb{K}[X_1,\dots,X_n]/I':
\end{equation*}
\begin{enumerate}
\item $\phi$ is an isomorphism;
\vskip2mm
\item there is a change of variables $F$ in $\mathbb{K}[X_1,\dots,X_n]$ such that $I'= F(I) + \langle X_1,\dots,X_n\rangle^l$ and $\phi(f+ I) = F(f) + I'$ for every $f\in \mathbb{K}[X_1,\dots,X_n]$.  
\end{enumerate}
In particular, when $I = I'$, every automorphism of A is induced by a change of variables F such that $F(I)\subset I$.
\end{proposition}
\vskip2mm
\noindent
We will denote by $C_{(A, I)}$ (resp. $CL_{(A, I)}$), the set of changes of variables (resp. linear changes of variables) $F$ of $\mathbb{K}[X_1,\dots,X_n]$ satisfying $F(I)\subset I$. One knows that $C_{(A,I)}$ is a monoid and $CL_{(A,I)}$ is isomorphic to a subgroup $L_A$ of $G_A$ such that $L_A\cap \text{Ker}(\Phi_A)=\{1\}$ (see \cite{AS3}, Corollary 21).
\vskip2mm
\noindent
Next, we review Lemma $22$ of \cite{AS3} in the context of local commutative algebras $A$ over $\mathbb{K}$.
\begin{lemma}\label{Lemma 4.3}
$($\cite{AS3}, \emph{Lemma} $22)$ Let $A$ be a local commutative algebra over $\mathbb{K}$ with $\emph{dim}(R/R^2)=n$. Suppose $A\cong \mathbb{K}[X_1,\dots,X_n]/I$ be a quiver representation with $I$ an admissible ideal (a finite set of adequate relations). If $I$ is a homogeneous ideal (i.e., $A$ is graded by the radical algebra), i.e, if homogeneous elements generate the admissible ideal $I$, then
\begin{equation*}\label{eq20}
\emph{\text{Im}}(\Phi_A) \cong CL_{(A, I)}\cong L_A.
\end{equation*}
\end{lemma}
\vskip2mm
\noindent
Now we mention a result from \cite{DD} which is under preparation. For the sake of completeness, we include the proof.
\vskip2mm
 \begin{proposition}\label{Proposition 4.4}$($\cite{DD}, \emph{Proposition} $5.4)$
Let $A$ be a local commutative algebra that is graded by the radical over $\mathbb{K}$ with $\emph{dim}(R/R^2)=n$. Then, \[\emph{\text{Aut}}_{\mathbb{K}}(A)\cong \emph{\text{Im}}(\Phi_A)\ltimes \emph{\text{Ker}}(\Phi_A), \text{ where $\Phi_A$ is the canonical map as in Proposition \ref{Proposition 3.12}. }\]
\end{proposition}
\vskip1mm
\begin{proof}
We have the canonical homomorphism $\Phi_A: \text{Aut}_{\mathbb{K}}(A)\rightarrow \text{GL}(R/R^2)$ whose kernel is unipotent for local algebra $A$ by Proposition \ref{Proposition 3.12}. Therefore, we have the following short exact sequence: 
\begin{center}
\begin{tikzcd}
1\arrow[r] & \text{Ker}(\Phi_A)\arrow[r] & \text{Aut}_{\mathbb{K}}(A)\arrow[r] & \text{Im}(\Phi_A)\arrow[bend left=33]{l}{s}\arrow[r] & 1
\end{tikzcd}
\end{center}
It is known from Lemma \ref{Lemma 4.3} that any element of $\text{Im}(\Phi_A)$ gives an automorphism of $A$, which is induced by a linear change of variables. Hence, we have a section $s: \text{Im}(\Phi_A)\rightarrow \text{Aut}_{\mathbb{K}}(A)$ given by $s(M)=\phi_M$, where $\phi_M(f+I)=f(XM)+I$, where $X=(X_1,\dots,X_n)$. Finally, we have $\Phi_A\circ s=Id$, which further implies that the above sequence is split exact and provides the required isomorphism.
\end{proof}
\vskip2mm
\noindent
Therefore, for graded by the radical local algebra $A$, $G_A$ is solvable if and only if $\text{Im}(\Phi_A)$ is solvable.
\vspace{0.1in}
\begin{proposition}\label{Proposition 4.5}
Let $A$ be a local commutative graded by the radical algebra over $\mathbb{K}$ with $\emph{dim}(R/R^2)=n$. Then there is a central torus $\mathbb{G}_m$ in $\emph{\text{Im}}(\Phi_A)$. It embeds as the diagonal change of variables in every maximal torus of $G_A$, given by \[\{\emph{diag }(a,\dots, a):a\in\mathbb{G}_m\}.\]
\end{proposition}
\vskip2mm
\begin{proof}
Let $A$ be graded by the radical algebra, which implies that there exists a homogeneous admissible ideal $I$ of $A$. We know that any automorphism $\phi$ of $A$ is induced by a change of variables $F$ such that $\phi(f+I)=F(f)+I$ and $F(I)\subseteq I$. Consider the following embedding $\Gamma :\mathbb{G}_m\rightarrow G_A$ given by $\Gamma(\alpha)=\phi_{\alpha}$, where $\phi_{\alpha}(\bar X_i)=\alpha \bar X_i$, where $\bar X_i$ denotes $X_i+I$. Therefore, for each $\phi_{\alpha}$, the corresponding change of variables is $F_{\alpha}(X_i)=\alpha X_i$, which stabilizes $I$. This is an injective homomorphism. If $\alpha\in \mathbb{G}_m$ is in $\text{Ker}(\Gamma)$, then $(\alpha-1)X_i\in I\subseteq \langle X_1,\dots,X_n\rangle^2, \forall i$, which is only possible if $\alpha=1$. Moreover, $\mathbb{G}_m$ embeds in $\text{Im}(\Phi_A)$ as $\text{Ker}(\Phi_A)$ is unipotent. 
\vskip2mm
\noindent
Let $\phi_1\in G_A$ be such that it is given by a linear change of variables, i.e., $\phi_1(f+I)=F_1(f)+I$, where $F_1$ is a linear change of variables corresponding to $\phi_1$, and let $\phi_2\in \mathbb{G}_m$, with $F_2(X_i)=\alpha X_i$ as the corresponding change of variables. Then, it is easy to verify that $F_1\circ F_2=F_2\circ F_1$, as $F_1(f)(X)=f(\sum_{i=1}^n\lambda_{i1}X_i,\dots,\sum_{i=1}^{n}\lambda_{in}X_i)$, where $(\lambda_{ij})_{n\times n}$ is a non-singular matrix. From Lemma \ref{Lemma 4.3} it follows that $\text{Im}(\Phi_A)$ is given by a linear change of variables. Hence, via the above embedding, $\mathbb{G}_m$ is a central torus in $\text{Im}(\Phi_A$). However, for any maximal torus $T$ in $G_A$, $T.\mathbb{G}_m$ is again a torus in $G_A$ containing $T$. Therefore, $\mathbb{G}_m\subset T$.  
\end{proof}

\vskip2mm
\begin{lemma}\label{Lemma 4.6}
Let $A$ be a local commutative graded by the radical algebra over $\mathbb{K}$. Suppose $S$ and $T$ are subgroups of $\emph{\text{Im}}(\Phi_A)$. Then, the central $\mathbb{G}_m$ (defined in Proposition \ref{Proposition 4.5}) of $\emph{\text{Im}}(\Phi_A)$ is not contained in $[S, T]^n$,  for any $n\in \mathbb{N}$. In particular, $[\emph{\text{Im}}(\Phi_A), \emph{\text{Im}}(\Phi_A)]\neq \emph{\text{Im}}(\Phi_A)$.
\end{lemma}
\vskip2mm
\begin{proof}
It is easy to see that $\text{det}([S,T]^n)=1$, where $\text{det}$ is the usual determinant map and $n\in\mathbb{N}$. If $\mathbb{G}_m\subset [S,T]^n$, then for each $\alpha\in \mathbb{K}^{\times}$, we have $\alpha^s=1$, where $s=\text{dim}( R/R^2)$, which is a absurd. Since $\mathbb{G}_m\subset \text{Im}(\Phi_A$), but $\mathbb{G}_m \nsubseteq [\text{Im}(\Phi_A), \text{Im}(\Phi_A)]$, we have $[\text{Im}(\Phi_A),\text{Im}(\Phi_A)]\neq \text{Im}(\Phi_A)$.
\end{proof}
\vskip2mm
\noindent
Again, we recall Lemma $5.9$ from \cite{DD} and include the proof for completeness.
\vskip2mm
\begin{lemma}\label{Lemma 4.7}$($\cite{DD}, \emph{Lemma} $5.9)$
Let $A$ be a local commutative algebra over $\mathbb{K}$ with $\emph{dim}(R/R^2)=n$, and by quiver and relations representations we have an isomorphism $A\cong \mathbb{K}[X_1,\dots,X_n]/I$, where $\langle X_1,\dots,X_n \rangle^l \subseteq I\subseteq \langle X_1,\dots,X_n \rangle^2$, $l\geq 2$, $I$ is an admissible ideal of $A$. Then we can write $I=\langle X_1,\dots, X_n \rangle^l+\langle P_1,\dots, P_m \rangle$, where $l$ is the Lowey length of $A$, and all non-zero $P_{i}$'s are polynomials of degree at least $2$ and at most $l-1$ with no linear component.
\end{lemma}
\vskip1mm
\begin{proof}
Let $f\in I$. Then $f=\sum_{i_1,\dots,i_n}a_{i_{1},\dots,i_{n}}X_{1}^{i_{1}}\dots X_{n}^{i_{n}}$. By the given hypothesis, we have $\langle X_1,\dots, X_n \rangle^l\subseteq I$; therefore, monomials with degree greater than or equal to $l$ are contained in $I$. Hence, we can write $f=\sum_{i}\lambda_{i}h_i+P$, where $h_i$'s are distinct monomials in \{$X_1,\dots,X_n$\} of degree at least $l$ and $P$ is a polynomial of degree at most $l-1$. However, $f-\sum_{i}\lambda_{i}h_{i}\in I$, therefore $P\in I$. Since $I$ is contained in $\langle X_1,\dots,X_n\rangle ^2$, $P$ must have degree at least $2$ with no linear component. Now, for each $f\in I$, we have a unique such $P$ in $I$ which has the degree at least $2$ and at most $l-1$ with no linear component. Consider the set
\[S=\{P: P \text{ corresponds to each } f\in I \text{ with the above property }\}.\] Suppose $I$ is a proper ideal lying in between $\langle X_1,\dots,X_n\rangle^l$ and $\langle X_1,\dots,X_n\rangle ^2$. Then $S$ must be non-empty. Therefore, the ideal generated by $S$ is contained in $I$, and let us call it $I_1$. Moreover, $k[X_1,\dots,X_n]$ is a noetherian ring, so $I_1= \langle P_1,\dots,P_m \rangle$ where $P_i\in S$ and $I_1\subset I$. Hence $I=\langle X_1,\dots,X_n\rangle^l+\langle P_1,\dots,P_m\rangle$, where ${P_i}$'s have the property given in the hypothesis.
\end{proof}

\vskip2mm
\begin{definition}\label{Definition 4.8} 
Stabilizer subgroup of a polynomial $f\in \mathbb{K}[X_1,\dots,X_n]$ in $\text{GL}(R/R^2)$, where $\text{dim}(R/R^2)=n$, is defined by 
\begin{equation*}\label{eq28}
\text{Stab}(f)\coloneq\{M\in \text{GL}(R/R^2): f(XM)=f(X), X=(X_1,\dots,X_n)\}.
\end{equation*}
\end{definition}
\noindent
$\text{Stab}^{0}(f)$ denotes the identity component of the stabilizer subgroup in $\text{GL}(R/R^2)$.

\vskip6mm
\begin{definition}\label{Definition 4.9}
The similarity subgroup of a polynomial $f\in \mathbb{K}[X_1,\dots,X_n]$ in $\text{GL}(R/R^2)$, where $\text{dim}(R/R^2)=n$, is defined by \[\text{Sim}(f)\coloneq\{M\in \text{GL}(R/R^2): f(XM)=\alpha_M f(X), X=(X_1,\dots,X_n)\text{ and } \alpha_M\in \mathbb{K}^{\times}\}.\]
\end{definition}
\noindent
$\text{Sim}^{0}(f)$ denotes the identity component of the similarity subgroup in $\text{GL}(R/R^2)$. It is easy to see that $\text{Sim}^{0}(f)=\mathbb{G}_m.\text{Stab}^{0}(f)$.

\vskip2mm
\begin{proposition}\label{Proposition 4.10}
Let $A$ be a local commutative graded by the radical algebra with $\emph{dim}(R/R^2)=n$. Suppose the corresponding admissible ideal for $A$ is $I\cong \langle X_1,\dots, X_n \rangle^l+\langle f \rangle$, where $2\leq \emph{deg }(f)<l$ is the Lowey length of $A$ and $f$ is a homogeneous polynomial.  Then,
\begin{equation*}\label{eq23}
\emph{\text{Aut}}_{\mathbb{K}}(A)\cong \emph{Sim}(f)\ltimes \emph{\text{Ker}}(\Phi_A), \text{ 
where $\Phi_A$ is the canonical map as in Proposition \ref{Proposition 3.12} }.
\end{equation*}
Therefore, $G_A$ is solvable if and only if $\emph{Stab}(f)^0$ is solvable.
\end{proposition}
\vskip1mm
\begin{proof}
We have the canonical map $\Phi_A: \text{Aut}_{\mathbb{K}}(A)\rightarrow \text{GL}(R/R^2)$. It is enough to prove that for the given $I$, $\text{Im}(\Phi_A)=\text{Sim}(f)$. Let $M\in \text{Im}(\Phi_A)$, then $M$ induces an automorphism of  $A$ which is given by the linear change of variables $F_1(X)=XM$, where $X=(X_1,\dots,X_n)$ and $F_1(I)\subset I$ by Lemma \ref{Lemma 4.3}. As $f\in I$,\[F_1(f)(X)=f(\sum_{i}m_{i1}X_i,\dots,\sum_{i}m_{in}X_i)\in I,\] where $M=(m_{ij})_{n\times n}$ and $X=(X_1,\dots,X_n)$. Since a linear change of variables does not change the degree of the polynomial, we have $\text{deg}(F_1(f))=\text{deg}(f)$ and $F_1(f)$ is a homogeneous polynomial. Therefore, $F_1(f)=gf$, where $g\in \mathbb{K}[X_1,\dots,X_n]$ is a non-zero homogeneous polynomial as $f$ is a homogeneous polynomial. However, $\text{deg}(gf)=\text{deg}(g)+\text{deg}(f)=\text{deg}(f)$ implies $\text{deg}(g)=0.$
Therefore, $F_1(f)=\alpha f$ and $f(\sum_{i}m_{i1}X_i,\dots,\sum_{i}m_{in}X_i)=\alpha f(X_1,\dots,X_n)$ with $\alpha\in \mathbb{K}^{\times}$.
Hence, $M\in \text{Sim}(f)$, so we have $\text{Im}(\Phi_A)\subset \text{Sim}(f)$. It is easy to see that $\text{Sim}(f)\subset \text{Im}(\Phi_A)$. Now using Proposition \ref{Proposition 4.4}, we are done. Hence, we have $G_A\cong \mathbb{G}_m.\text{Stab}^{0}(f)\ltimes \text{Ker}(\Phi_A)$, which implies $G_A$ is solvable if and only if $\text{Stab}^{0}(f)$ is solvable.
\end{proof}
\vskip2mm
\noindent
So, to study the automorphism group of local commutative associative algebras, we see that it is important to examine the stabilizer of polynomials where these polynomials occur as the generator of the ideal $I$ of $A$ in the quiver and relation representation of $A$. Therefore, we now mention an important result related to the stabilizer of a polynomial from the paper \cite{GG}.

\vskip2mm
\begin{definition}\label{Definition 4.11} 
Stabilizer of a polynomial $f\in \mathbb{K}[V]$ is defined by 
\begin{equation*}\label{eq28}
\text{Stab}(f)=\{g\in \text{GL}(V): f\circ g=f\}
\end{equation*}
\end{definition}

\vskip2mm
\begin{lemma}\label{Lemma 4.12}
$($\cite{GG}, \emph{Lemma} $5.1)$ Let $X\subset \emph{SL}(V)$ be a simple algebraic group over an algebraically closed field $\mathbb{K}$ such that $V$ is irreducible, restricted, and tensor indecomposable for $X$. Put $q$ for a non-zero $X$ invariant quadratic form on $V$ if one exists; otherwise, set $q\coloneq0$. If $(X,V)$ does not appear in Table $1$ of \cite{Sei}, then, for every $f\in (\mathbb{K}[V]^{X}-\mathbb{K}[q])$, the stabilizer of $f$ in $\emph{GL}(V)$ has identity component $X$.
If additionally $\emph{\text{char }}\mathbb{K}=p\neq 2,3$ and does not divide $\emph{deg}(f)$, $f$ is not in $\mathbb{K}[V]^{(p)}[q]$ (see \cite{GG} for notation), and furthermore $\emph{\text{char }}\mathbb{K}$ does not divide $n+1$ if $X$ has type $\emph{A}_n$, then the scheme-theoretic stabilizer of $f$ in $\emph{GL}(V)$ is smooth with identity component $X$.
\end{lemma}

\vskip2mm
\begin{remark}\label{Remark 4.13}
Over a field of characteristic $0$, any irreducible representation is restricted and tensor indecomposable (see \cite{GG}, \cite{Sei}). Therefore, we can construct different examples of non-solvable groups $G_A$ using the Lemma \ref{Lemma 4.12} and Proposition \ref{Proposition 4.10} over an algebraically closed field with characteristic $0$. 
\end{remark}
\vskip2mm
\begin{definition}\label{Definition 4.14}
$($\cite{OS}, \cite{Schneider}$)$
    A non-zero homogeneous polynomial $f\in \mathbb{K}[X_1,\dots, X_n]$ of degree $d$ is called non-singular if any one of the following equivalent conditions holds;
    \begin{enumerate}
        \item $0=(0,\dots,0)$ is the unique common root of its partial derivatives $\displaystyle\frac{\partial f}{\partial X_1},\dots,\displaystyle\frac{\partial f}{\partial X_n}$. 
        \vskip2mm
        \item Discriminant $(\Delta_{n,d}(f))$ of $f$ does not vanish (see Chapter $13$, \cite{GKZ}).
        \vskip2mm
        \item Let $\Theta_f$ be the corresponding symmetric multilinear form of $f$ of degree $d$ over a vector space $V=\mathbb{K}^n$ of dimension $n$. Then, $\Theta_f(v,\dots,v,w)=0$ for all $w\in V$ implies $v=0$. 
    \end{enumerate}
\end{definition}
\vskip2mm
\noindent
Next, we recall the following theorem, originally due to Camille Jordan.
\vskip2mm
\begin{theorem}\label{Theorem 4.15}
$($\cite{Schneider}$)$ Let $\Theta$ denote a symmetric multilinear map of degree $r\geq 3$ on a vector space of dimension $n$ over a field $K$ (arbitrary field). Assume that $K$ has the characteristic $0$ or greater than $r$. If $\Theta$ is non-singular, then its orthogonal group (stabilizer group) in $\emph{GL}_n(K)$ is finite.
\end{theorem}
\vskip2mm
\noindent
We now prove the following result:
\vskip2mm
\begin{proposition}\label{Proposition 4.16}
Let $A$ be a local commutative graded by the radical algebra over $\mathbb{K}$ with Lowey length $l$ with $\emph{dim}(R/R^2)=n$, where the admissible ideal $I$ is given by $\langle X_1,\dots, X_n \rangle^l +\langle f_1,\dots,f_m\rangle$, $f_i$'s are homogeneous polynomials such that not all $f_i$'s are zero. Let $d_i=\emph{\text{deg}}(f_i)$, where $1\leq i\leq m$ and $2\leq d_i<l$.  Suppose there exists a unique generator $f_s$ with degree $\emph{\text{min}}\{d_i: 1\leq i\leq m\}$. Then 
\begin{equation*}\label{eq36}
\bigcap_{i=1}^{m}\emph{Sim}(f_i)\subset\emph{\text{Im}}(\Phi_A)\subset \emph{Sim}(f_s). 
\end{equation*} 
Hence, if $G_A$ is solvable, then $\bigcap_{i=1}^{m} \emph{Sim}(f_i)^{0}$ is solvable. If $\emph{Sim}(f_s)$ is solvable, then $G_A$ is solvable.
In particular, if $f_s$ is a non-singular homogeneous polynomial of degree $\geq 3$ over $\mathbb{K}$, where $\emph{char }\mathbb{K}=0$ or $\emph{char }\mathbb{K}=p>3$ and $f_s$ is a non-singular homogeneous polynomial of degree $d$ with $2<d<p$, then $G_A$ is a solvable group of rank $1$.
\end{proposition}
\vskip1mm
\begin{proof}
 It is known from Lemma \ref{Lemma 4.3} that for a graded by the radical algebra $A$, \[\text{Im}(\Phi_A)=\{M\in \text{GL}(R/R^2): \phi_M\in \text{Aut}_{\mathbb{K}}(A) \text{ with } \phi_M(f+I)=f(XM)+I \}\] Therefore, every $M\in \text{Im}(\Phi_A)$ induces an automorphism of $A$ which is given by the linear change of variables $F_1(X)=XM$, where $X=(X_1,\dots,X_n)$ and $F_1(I)\subset I$ by Proposition \ref{Proposition 4.2} and Lemma \ref{Lemma 4.3}. It is easy to see that $\bigcap_{i=1}^{m}\text{Sim}(f_i)\subset \text { Im}(\Phi_A)$. Let us denote $f_M(X)\coloneq f(XM)$ and a linear change of variables by  $F_1$ for $M$, where $M\in \text{ Im}(\Phi_A)$. Since a linear change of variables does not change the degree of a polynomial and $F_1(f_i)={f_{i}}_{M} \in I$ by Lemma \ref{Lemma 4.3}, it follows that ${f_{i}}_M$ is a homogeneous polynomial of degree $d_i$ and ${f_{i}}_M=F_1(f_i)=\sum_{j=1}^{m}g_{ij}f_j$, where $g_{ij}\in \mathbb{K}[X_{1},\dots,X_{n}]$ is a homogeneous polynomial of degree $s_{ij}$. As $F_1$ is given by the linear change of variables, we have
\begin{equation*}\label{eq30} 
d_i= \text{deg}({f_{i}}_{M})=\text{deg}(F_1(f_i))=\text{deg }(\sum_{j=1}^{m}g_{ij}f_j)=\text{ Max}_{j=1}^{m}(\text{deg}(g_{ij})+\text{deg}(f_j)).
 \end{equation*}
 Without loss of generality, let $d_1=\text{ min }d_i$, and $d_{i}\leq d_{i+1}$ for all $i\geq 2$, where $d_i=\text{deg}(f_i)$. Since ${f_{i}}_M$ is a homogeneous polynomial of degree $d_i$, we have $d_i=(s_{ij}+d_j)$ for all $1\leq i,j \leq m$, where $\text{deg}(g_{ij})=s_{ij}$. By the given hypothesis, $f_1$ is the unique generator of minimum degree $d_1$. Therefore, $d_1=s_{11}+d_1$ and $d_1<d_i$ for all $i\neq 1$, which implies $s_{11}=0$ and $g_{1j}=0,\forall j\in\{2,\dots,m\}$. So, we get $F_1(f_1)(X)=f_1(XM)=\alpha f_1$. Hence, $\text{Im}(\Phi_A)\subset \text{Sim}(f_1)$. Finally, we have \[\bigcap_{i=1}^{m}\text{Stab}(f_i)\subset\bigcap_{i=1}^{m} \text{Sim}(f_i)\subset \text{Im}(\Phi_A)\subset \text{Sim}(f_1).\] However, $\text{Sim}^{0}(f)=\mathbb{G}_m.\text{Stab}^{0}(f)$. Therefore, we get a necessary condition for the solvability of $G_A$, that is, $G_A$ is solvable implies that $\bigcap_{i=1}^{m} \text{Stab}(f_i)^0$ is solvable. If $\text{Sim}(f_1)$ is solvable, then $\text{Im}\Phi_A$ is solvable, which implies $G_A$ is solvable.
\vskip2mm
\noindent
It follows from Proposition \ref{Proposition 4.4} that $G_A\cong \text{Im}(\Phi_A)^{0}\ltimes\text{Ker}(\Phi_A)^{0}$, and $\text{Ker}(\Phi_A)^{0}$ is connected unipotent group.  It is clear from Theorem \ref{Theorem 4.15} that under the assumption on $f_1$ depending of characteristic of $k$, $\text{Sim}(f_1)^{0}=\mathbb{G}_m$ as $\text{Stab}^{0}(f_1)=\{Id\}$. Hence, under the assumption on $f_1$, $G_A$ is a solvable group of rank $1$. 
\end{proof}

\vskip2mm
\noindent
Next, we recall the following important lemma from \cite{DD} and include the proof for completeness.
\vskip2mm
\begin{lemma}\label{Lemma 4.17}$($\cite{DD}, \emph{Lemma} $5.30)$
Let $A$ be a local commutative graded by the radical algebra over $\mathbb{K}$ with Lowey length $l$ and $\emph{dim}(R/R^2)=n$. Then the admissible ideal $I$ contains an $\emph{\text{Im}}(\Phi_A)$-stable finite-dimensional subspace $W$ with all non-zero elements having the same degree.
\end{lemma}
\vskip1mm
\begin{proof}
Since $A$ is a graded by the radical local algebra, $I$ is isomorphic to $\langle X_1,\dots,X_n\rangle^l+\langle f_1,\dots,f_m\rangle$ with $f_i$' s being homogeneous polynomials of $2\leq\text{deg}(f_i)\leq (l-1)$. Let $S$ be the collection of polynomials of minimal degree in the generating set of $I$; \[S=\{f_i: \text{deg}(f_i)\leq \text{deg}(f_j), \forall i\neq j \text{ and }i,j\in\{1,\dots,m\}\}.\] By renumbering if necessary, we can assume that $S=\{f_1,\dots,f_r\}$ for some $r\leq m$. We claim that $W\coloneq\text{ Span }S$ is an $\text{Im}(\Phi_A)$-stable subspace under the action defined by $F_1.f=F_{1}(f)$, for all $F_1\in \text{Im}(\Phi_A)$ and $f\in\mathbb{K}[X_1,\dots,X_n]$, where $F_{1}(f)(X)=f(XM_{F_{1}})$. We have $W=\{c_1f_1+\dots+c_rf_r:c_i\in \mathbb{K}\}$. It is known that $F_1(I)\subset I$ for all $F_1\in \text{Im}(\Phi_A$) by Lemma \ref{Lemma 4.3}, where $F_1$ is given by a linear change of variables. Let $0\neq g\in W$, then $F_1(g)\in I$. We have $\text{deg}(F_{1}(g))=\text{deg}(g)$ and $F_{1}(g)$ is a non-zero homogeneous polynomial as $F_1$ is given by a linear change of variables. Hence, $\text{deg}(g)=\text{deg}(F_1(g))=\text{Max}_{j=1}^m(\text{deg}(g_j)+\text{deg}(f_j))=\text{deg}(g_j)+\text{deg}(f_j)$, where $F_1(g)=\sum_{j=1}^{m}g_jf_j$, $g_j\text{ is a homogeneous polynomial in }\mathbb{K}[X_1,\dots,X_n]$ and $j\in\{1,\dots,m\}$. However, $\text{deg}(g)=\text{deg}(f_i)<\text{deg}(f_s)$, $\forall i\in \{1,\dots,r\},\forall s\in \{r+1,\dots,m\}$, therefore, $\text{deg}(g_i)=0$ and $g_s=0$. This implies $F_1(g)=a_1f_1+a_2f_2+\dots+a_rf_r\in W$, where $a_i\in \mathbb{K}.$
\end{proof}
\vskip2mm
\begin{remark}\label{Remark 4.18}
    The above subspace $W$ is independent of the choice of generators of $I$. We call this subspace the \emph{minimal degree subspace} associated to $A$ as it is generated by polynomials of minimal degree in $I$ and $I$ is determined by $A$. We denote the \emph{minimal degree subspace} associated to $A$ by $W$.
\end{remark}
\vskip2mm
\noindent
We need the following result on closed orbits from \cite{Humphreys}.
\vskip2mm
\begin{proposition}\label{Proposition 4.19}
$($\cite{Humphreys}, \emph{Proposition} $8.3)$
Let the algebraic group $G$ act morphically on the non-empty variety $X$. Then each orbit is a smooth, locally closed subset of $X$, whose boundary is the union of orbits of strictly lower dimension. In particular, orbits of minimal dimension are closed. (So, closed orbits exist). 
\end{proposition}
\vskip2mm
\noindent
Next, we recall the following result from \cite{kraft}.
\vskip2mm
\begin{lemma}$($\cite{kraft}, \emph{Corollary} $4.2)$\label{Lemma 4.20} Let $G$ be a connected reductive group and $\mathbb{K}[V]$ be the ring of regular functions on a $G$-module $V$ $(\emph{dim }V=n)$ and let $f_1,\dots,f_n$ be a regular sequence of homogeneous elements of $\mathbb{K}[V]$ such that the linear span  $\emph{Span }\{f_1,\dots,f_n\}$ is $G$-stable. Then, $\emph{Span }\{f_1,\dots,f_n\}$ contains all non trivial representations of $[G,G]$ in $\mathbb{K}[V]_1$, where $\mathbb{K}[V]$ has a natural grading by degree, i.e, $\mathbb{K}[V]=\bigoplus_{i} \mathbb{K}[V]_{i}$ and $\mathbb{K}[V]_{1}$ is the linear part of $\mathbb{K}[V]$.
\end{lemma}
\vskip2mm
\noindent
Let us recall that if $G$ be a connected algebraic group over an algebraically closed field $\mathbb{K}$ with $\text{dim }G\leq 2$, then $G$ is solvable (see \cite{Borel}, Corollary $11.6$). Now we prove the following lemma, which is used to prove the main result.
\vskip2mm
\begin{lemma}\label{Lemma 4.21}
Let $G$ be a connected affine algebraic group over an algebraically closed field $\mathbb{K}$. Then $G$ is solvable if and only if it does not contain any connected non-trivial semisimple subgroup.
\end{lemma} 
\vskip1mm
\begin{proof}
If $G$ is solvable, then it can not contain any semisimple subgroup, as a subgroup of a solvable group is solvable, but semisimple groups are not solvable. Conversely, we can show that if $G$ does not contain any connected non-trivial semisimple subgroup, then $G$ is solvable. Therefore, it is enough to show that any non-solvable connected group $G$ must contain a non-trivial connected semisimple subgroup. To prove this, we will use strong induction on the dimension of $G$. It is known that any connected group of dimension at most $2$ is solvable (see \cite{Borel}, Corollary $11.6$). Let us consider a group $G$ which is non-solvable and $\text{dim }G=3$ with $U$ its unipotent radical. We claim that $U=\{e\}$. Suppose $U\neq \{e\}$, then $\text{dim }U\geq 1$. It implies $\text{dim }G/U\leq 2$. Hence, $G/U$ is solvable. Therefore, $G$ is solvable as $U$ is solvable, which is a contradiction as $G$ is non-solvable. So we have $G$ is reductive, if $\text{dim }G=3$ and $G$ is non-solvable. However, from the structure theory of a reductive group, it is known that $[G, G]$ is a semisimple subgroup of $G$. Since $G$ is non-solvable, $[G,G]\neq \{e\}$. We claim that $G=[G,G]$. Suppose $G\neq [G,G]$, then $\text{dim }[G,G]\leq 2$ and $G/[G,G]$ is solvable. It implies $G$ is solvable, which is again a contradiction. Therefore, we have $G$ is semisimple, if $\text{dim }G=3$ and $G$ is non-solvable. So the strong induction process starts at dimension $3$. Let us assume $G$ is a connected non-solvable group of dimension $n\geq 3$. Suppose for any connected group $H$ whose dimension $d$, where  $3\leq d\leq (n-1)$ and $H$ is non-solvable, $H$ must contain a non-trivial connected semisimple subgroup. Now we prove that $G$ contains a non-trivial connected semisimple subgroup. Let $U$ be the unipotent radical of $G$. If $U=\{e\}$, then $G$ is reductive. Hence, by the structure theory of reductive groups, $G$ must have a non-trivial connected semisimple subgroup. So we may assume $U\neq \{e\}$, then consider $G/U$. Therefore, $3\leq \text{dim }G/U<\text{dim }G$ and $G/U$ is not solvable as $G$ is not solvable. Hence, by the strong induction hypothesis, we have that $G/U$ must contain a non-trivial connected semisimple subgroup $T$. However, we have the canonical surjective map $\pi: G\rightarrow G/U$. So $\pi^{-1}(T)$ is non-solvable as $T$ is semisimple and $3\leq \text{dim }\pi^{-1}(T)<\text{dim }G$. Therefore, by strong induction, $\pi^{-1}(T)$ must contain a non-trivial connected semisimple subgroup. Hence, $G$ contains a non-trivial connected semisimple subgroup.
\end{proof}

\vskip2mm
\noindent
From now onward, we assume that $\text{Im}(\Phi_A)$ is connected and for $f\in \mathbb{K}[X_1,\dots,X_n]$, the orbit of $f$ under the action of $\text{Im}(\Phi_A)$ is denoted by $O(f)\coloneq \{F.f: F\in \text{Im}(\Phi_A)\}$, where $(F.f)=F(f)$, and $F(f)(X)=f(XM_{F}), X=(X_1,X_2,\dots,X_n)$, $M_{F}$ is the corresponding linear change of variables matrix.
\begin{proposition}\label{Proposition 4.22}
  Let $A$ be a local commutative graded by the radical algebra over $\mathbb{K}$ with $\emph{dim}(R/R^2)=n$. Assume that $W$ is the minimal degree subspace associated to $A$. Suppose $\emph{\text{char }}\mathbb{K}=0$ and $W$ contains a non-singular homogeneous polynomial of degree at least $3$, additionally if $\emph{\text{char }}\mathbb{K}=p >3$ and $W$ contains a non-singular homogeneous polynomial of degree $d$ with $2<d<p$. Then $\emph{\text{dim }}\emph{\text{Im}}(\Phi_A)=\emph{\text{dim }}W$ and all non-singular elements in $W$ form an open $\emph{\text{Im}}(\Phi_A)$-orbit in $W$. Hence, $(\emph{\text{Im}}(\Phi_A),W)$ becomes a prehomogeneous vector space.
\end{proposition}
\vskip1mm
\begin{proof}
It is clear from Lemma \ref{Lemma 4.17} that, under the action of $\text{ Im}(\Phi_A)$, $W$ is an $\text{Im}(\Phi_A)$-stable subspace. Assume $f$ is a non-singular polynomial of degree $d\geq 3$ in $W$ if $\text{char }\mathbb{K}=0$ and if $\text{char }\mathbb{K}=p>3$, then $f$ is a non-singular polynomial of degree $d$ with $2<d<p$. Consider \[O(f)=\{F.f: F \in \text{Im}(\Phi_A)\}\] We have a bijection of $\text{Im}(\Phi_A)/\text{Stab}(f)\cong O(f)$ and 
\begin{equation}\label{eq 3.1.11}
\text{dim }\text{Im}(\Phi_A)-\text{dim }\text{Stab}(f)=\text{dim }O(f).
\end{equation}
\noindent
Moreover, by Theorem \ref{Theorem 4.15}, $\text{dim}(\text{Stab}(f))=0$ under the hypothesis on $W$ depending on the characteristic of $\mathbb{K}$, which implies that $\text{dim }O(f)=\text{dim }\text{Im}(\Phi_A$). Since $O(f)\subset W$, $\text{dim } O(f)\leq \text{dim } W$. Hence, we have $\text{dim }\text{Im}(\Phi_A)\leq \text{dim } W$.
\vskip2mm
\noindent
Let $U$ be the collection of all non-singular polynomials in $W$. Therefore, $U$ is an open set in $W$ as it is the complement of the zero set in $W$ given by the irreducible polynomial $\Delta_{n,d}$, where $\Delta_{n,d}$ is the discriminant (Chapter $13$, \cite{GKZ}). We claim that \[U=O(f) \text{ for some } f\in U.\] Let $f\in U$, then $\Delta_{n,d}(F.f)=\text{det}(M_{F})^t.\Delta_{n,d}(f)\neq0$ for some $t$, as $\text{det}(M_{F})\neq 0$ for all $F\in \text{Im}(\Phi_A)$ and $M_{F}$ is the corresponding linear change of variables matrix. Therefore, $U$ is stable under $\text{Im}(\Phi_A)$. Hence, $U=\bigsqcup_{f_{\alpha}\in U} O(f_{\alpha})$. However, $U$ is a noetherian topological space as it is a topological subspace of a noetherian space. Each $O(f_{\alpha})$ is irreducible as it is a continuous image of $\text{Im}(\Phi_A)$. Since each $f_{\alpha}$ is non-singular, it follows from Theorem \ref{Theorem 4.15} that $\text{dim}(\text{Stab}(f_{\alpha}))=0$. Therefore, using equation (\ref{eq 3.1.11}), we get $\text{dim }O(f_{\alpha})=\text{dim }O(f_{\beta})$ for all $f_\alpha$ and $f_\beta$ in $U$. Hence, each $O(f_{\alpha})$ is closed in $U$ because $\overline{O}(f_{\alpha})\cap U=O(f_{\alpha})$ by Proposition \ref{Proposition 4.19}, where $\overline{O(}f_\alpha)$ denotes the closure of $O(f_\alpha)$ in $W$. In a noetherian topological space, any nonempty closed subset can be written as a finite union of irreducible components. Therefore, $U=\bigsqcup_{i=1}^{n} O(f_i)$, for some $f_i\in U$, $1\leq i\leq n$. Since $U$ is an open subset of an irreducible variety $W$, $U$ is irreducible  as a topological space, hence $U=O(f_m)$ for some $m\in\{1,\dots,n\}$, and $W=\overline{U}=\overline{O}(f_m)$. Therefore, we have $\text{dim }\text{Im}(\Phi_A)=\text{dim }O(f_m)=\text{dim }U=\text{dim }W$.
\end{proof}
\vskip2mm
\begin{remark}\label{Remark 4.23}
Let $A$ satisfy the same hypothesis as in the above proposition. Then in this case, we can conclude that $\text{rank}(\text{Im}(\Phi_A))\leq \text{dim } W$. Therefore, $\text{rank} (G_A)\leq \text{dim }W$ as $\text{rank}(\text{Im}(\Phi_A))=\text{rank}(G_A)$ by Proposition \ref{Proposition 3.12}. Moreover, by Proposition \ref{Proposition 3.12}  we have $\text{rank}(G_A)\leq \text{dim}(R/R^2)$. However, $\text{rank}(G_A)\leq \text{dim }W$ is much better bound for the rank of $G_A$. For example, if $ A= \mathbb{K}[X_1,\dots,X_{20}]/I$ with $\text{char }\mathbb{K}=0$ or $\text{char }\mathbb{K}\geq 5$, where \[I=\langle X_1,\dots,X_{20}\rangle ^7+\langle X_1^4+X_2^4+\dots+X_{20}^4, X_1X_2X_3X_4, X_1^2X_5^2+X_2^2X_{17}^2\rangle,\] then $W=\text{Span }\{X_1^4+X_2^4+\dots+X_{20}^4, X_1X_2X_3X_4, X_1^2X_5^2+X_2^2X_{17}^2\}$. Hence, $\text{dim } W=3$ and $\text{dim}(R/R^2)=20$, therefore $\text{rank}(G_A)\leq 3$.
\end{remark}

\section{Proof of the Main Theorem}\label{sec: Main}
\noindent
The aim of this section is to prove the main theorem of this article and give some examples of new classes of algebras where we can apply the theorem. Now we prove the main result of this paper.
\vskip2mm
\begin{theorem}\label{Theorem 5.1}
Let $A$ be a local commutative graded by the radical algebra over an algebraically closed field $\mathbb{K}$ with $\emph{dim}(R/R^2)=n$. Then by quiver and relation representation of $A$, we have an isomorphism $A\cong \mathbb{K}[X_1,\dots, X_n]/I$  with $\langle X_1,\dots, X_n\rangle^l\subsetneqq I\subsetneqq \langle X_1,\dots, 
 X_n\rangle^2$, where $l>2$ is the Lowey length of $A$, and $I$ is a homogeneous ideal. Let $W$ be the minimal degree subspace associated to $A$, i.e., the collection of all non-zero homogeneous polynomials of minimal degree in $I$ and $0$ (see Remark \ref{Remark 4.18} ). Let any of the conditions listed below be satisfied:
    \vskip2mm
    \begin{enumerate}
    \vskip2mm
    \item $\emph{\text{char }}\mathbb{K}=0$ and $\emph{\text{dim }}W=n$ with $W$ having an ordered basis, whose elements form a regular sequence in the polynomial ring $\mathbb{K}[X_1,\dots, X_n]$;
    \vskip2mm
    \item $\emph{\text{char }}\mathbb{K}=0$ and $W$ contains a non-singular homogeneous polynomial of degree at least $3$;
\item $\emph{\text{char }}\mathbb{K}=p>3$ and $W$ contains a non-singular homogeneous polynomial of degree $d$ with $2<d<p$.  
    \end{enumerate}
    \vskip2mm
 Then, $G_A$ is solvable.
\end{theorem}
\vskip1mm
\begin{proof}
\vskip1mm
\noindent
\textbf{Proof of $(1)$:}
It is known that $G_A$ is solvable if and only if $\text{ Im}(\Phi_A)$ is solvable by Proposition \ref{Proposition 4.4}. It follows from Lemma \ref{Lemma 4.17} that $W$ is an $\text{Im}(\Phi_A)$-stable subspace. We prove the result by contradiction. Suppose $\text{ Im}(\Phi_A)$ is not solvable; then it must contain a non-trivial connected semisimple subgroup $H$ by Lemma \ref{Lemma 4.21}. Since $\text{Im}(\Phi_A)$ has an action on $\mathbb{K}[X_1,\dots,X_n]$, we can also consider $\mathbb{K}[X_1,\dots, X_n]$ as an $H$-module. Let $f_1,\dots,f_n$ be a regular sequence of homogeneous polynomials in $\mathbb{K}[X_1,\dots,X_n]$ that forms a basis in $W$. From Lemma \ref{Lemma 4.17}, it follows that $W$ is also an $H$-stable subspace of $\mathbb{K}[X_1,\dots,X_n]$. Therefore, using Lemma \ref{Lemma 4.20}, we can conclude that $W$ must contain all non-trivial representations of $H$ in $\mathbb{K}[X_1,\dots, X_n]_1$ (linear part), which is a contradiction as $\text{ deg}(f_i)\geq 2$.

\vspace{0.2in}
\noindent
\textbf{Proof of $(2)$ and $(3)$:}
If $\text{char }\mathbb{K}=0$, then $W$ contains a non-singular homogeneous polynomial of degree $d\geq 3$ and if $\text{char }\mathbb{K}=p>3$, then $W$ contains a non-singular homogeneous polynomial of degree $d$ with $2<d<p$. Therefore, we have $U=O(f)$ by Proposition \ref{Proposition 4.22}, where $U$ is the collection of non-singular elements in $W$ and $f$ is a non-singular homogeneous polynomial in $W$. Note that $U$ is an irreducible space as it is a non-empty open subset of $W$, and $W$ is a vector space (an irreducible affine variety). It is enough to prove $\text{Im}(\Phi_A)$ is solvable. If $\text{Im}(\Phi_A)$ is not solvable, then there exists a non-trivial connected semisimple subgroup $\tilde{G}$ of $\text{Im}(\Phi_A)$ by Lemma \ref{Lemma 4.21}. Let $f\in U$, then $\Delta_{n,d}(F_{1}.f)=\text{det}(M_{F_{1}})^m.\Delta_{n,d}(f)\neq0$ for some $m$, as $\text{det}(M_{F_{1}})\neq 0$ for all $F_{1}\in\tilde{G}$ and $M_{F_{1}}$ is corresponding linear change of variables matrix. Therefore, $U$ is stable under $\tilde{G}$. Consider the orbit decomposition of $U$ under the action of the subgroup $\tilde{G}$, hence we have $U=\bigsqcup_{f_\alpha\in U} O_{\tilde{G}}(f_\alpha)$. It is clear that $U$ is a noetherian, irreducible topological space. Moreover, by using similar equation like equation (\ref{eq 3.1.11}) for $\tilde{G}$, we get $\text{dim }(O_{\tilde{G}}(f_\alpha))=\text{dim }(O_{\tilde{G}}(f_\beta))$, which implies that $O_{\tilde{G}}(f_\alpha)$ is closed in $U$ by  Proposition \ref{Proposition 4.19}. Each $O_{\tilde{G}}(f_\alpha)$ is irreducible as it is a continuous image of $\tilde{G}$. Therefore, we have $U=\bigsqcup_{i=1}^{n}O_{\tilde{G}}(f_i)$ for some $f_i\in U$, $1\leq i\leq n$ as $U$ is a noetherian topological space. As a result, we have $U=O_{\tilde{G}}(f)$ for some $f\in U$, as $U$ is irreducible. However, we know that $\text{dim }\tilde{G} -\text{ dim }\text{Stab}_{\tilde{G}}(f)=\text{dim }O_{\tilde{G}}(f)$. Since $f$ is non-singular, using Theorem \ref{Theorem 4.15} we can conclude that $\text{dim }\tilde{G}=\text{dim }O_{\tilde{G}}(f)=\text{dim }U$ as $\text{dim }\text{Stab}_{\tilde{G}}(f)=0$ under the assumption on $W$. Again, $\text{dim Im}(\Phi_A)=\text{dim } U=\text{dim }{\tilde{G}}$. So, $\text{Im}(\Phi_A)=\tilde{G}$. Hence, $\text{Im}(\Phi_A)$ is semisimple which implies $[\text{Im}(\Phi_A),\text{Im}(\Phi_A)]=\text{Im}(\Phi_A)$, which is a clear contradiction by Lemma \ref{Lemma 4.6}. Hence, in both cases $\text{Im}(\Phi_A)$ has no connected non-trivial semisimple subgroup. It implies that $\text{Im}(\Phi_A)$ is solvable by Lemma \ref{Lemma 4.21}. 
\end{proof}
\vskip2mm
\begin{proposition}$($\cite{Popov}, \emph{Theorem} $10)$\label{Proposition 5.2}
    The following conditions are equivalent for a connected affine algebraic group $G$ over $\mathbb{C}$.
    \begin{enumerate}
        \item $G$ is solvable.
        \vskip2mm
        \item $m(G)=\emph{ rank }(Hom(G,G_m))$, where $m(G)=\emph{\text{max }}\{d\in \mathbb{N}\cup\{0\}:H_{d}(G,\mathbb{Q})\neq 0\}$, and $H_d(G,\mathbb{Q})$ is the singular homology considering $G$ as a complex manifold.
    \end{enumerate}
\end{proposition}
\vskip2mm
\begin{corollary}\label{Corollary 5.3}
Let $A$ be a graded by the radical local commutative algebra over $\mathbb{C}$. Then $G_A$ is solvable if and only if $m(G_A)=\emph{ rank }(G_A)$. In particular, if $W$ contains a non-singular element of degree at least $3$, where $W$ is defined as in Lemma \ref{Lemma 4.17}, then $H_d(G_A,\mathbb{Q})=0$, where $d>\emph{\text{min }}\{\emph{dim}(R/R^2),\emph{dim }W\}.$
\end{corollary}
\vskip2mm
\begin{proof}
 The first part of the corollary directly follows from Proposition \ref{Proposition 5.2}. By Proposition \ref{Proposition 3.12}, we know that $\text{ rank }(G_A)=\text{rank}(\text{Im}(\Phi_A))\leq \text{dim}(R/R^2)$. Let $W$ contain a non-singular homogeneous polynomial of degree at least $3$. Then, by Proposition \ref{Proposition 4.22}, we have  $\text{rank}(\text{Im}(\Phi_A)) \leq \text{dim }W$. Hence, $\text{rank}(G_A)\leq \text{dim }W$. Now, using Proposition \ref{Proposition 5.2} and Theorem \ref{Theorem 5.1}, we are done.
\end{proof}
\vskip2mm
\begin{example}\label{Example 5.4}
Consider the following example over the field $\mathbb{K}=\mathbb{C}$. Let \[A= \mathbb{C}[X_1,\dots,X_{25}]/I,\] where \[I=\langle X_1,\dots,X_{25}\rangle ^7+\langle X_1^3+X_1^2X_2+X_2^3+X_3^3+\dots+X_{25}^3, X_1X_2X_3, X_1X_5^2+X_2X_{17}^2\rangle,\] then $W=\text{Span }\{X_1^3+X_1^2X_2+X_2^3+X_3^3+\dots+X_{25}^3, X_1X_2X_3, X_1X_5^2+X_2X_{17}^2\}$. Hence, $\text{dim } W=3$, $\text{dim}(R/R^2)=25$ and $W$ contains a non-singular homogeneous polynomial of degree $3$, $f=X_1^3+X_1^2X_2+X_2^3+X_3^3+\dots+X_{25}^3$. By Proposition \ref{Proposition 3.12}, we have $\text{rank} (G_A)=\text{rank }(\mathrm{Im}(\Phi_A))$, and rank $(\mathrm{Im}(\Phi_A))\leq \text{dim }W$ by Proposition \ref{Proposition 4.22}. Therefore, we have  $\text{rank}(G_A)\leq \text{dim }W=3$. Using Theorem \ref{Theorem 5.1}, $G_A$ is solvable as $W$ contains a non-singular homogeneous polynomial $f=X_1^3+X_1^2X_2+X_2^3+X_3^3\dots+X_{25}^3$ of degree $3$. Now using Corollary \ref{Corollary 5.3}, we can conclude that $H_{d}(G_A,\mathbb{Q})=0$ for $d>3$.
\end{example}
\vskip2mm
\noindent 
We now present some examples of new classes of algebras over an algebraically closed field $\mathbb{K}$ that satisfy the conditions of the main theorem. These algebras are not included in the established results (Theorem \ref{Theorem 3.6}, Theorem \ref{Theorem 3.7}, Theorem \ref{Theorem 3.8}, Theorem \ref{Theorem 3.9}) previously mentioned in Section \ref{sec: Solvability}.
\vskip2mm
\noindent
\begin{example}\label{Example 5.5}
Let $A=\mathbb{K}[X_1,X_2,X_3]/I$, where \[I=\langle X_1, X_2, X_3\rangle^l+\langle X_1^2+X_2^2, X_2^2+X_3^2, X_2X_3, X_1X_2^2+X_3^3, X_1^4+X_2^4\rangle\]
with $l\geq 5$ and $\text{char }\mathbb{K}=0$. It is easy to see that $A$ is a local non-monomial graded by the radical algebra as $A/R=\mathbb{K}$, and not all the generators of the homogeneous ideal $I$ are monomials. Since $\mathbb{G}_m$ embeds in $G_A$ by Proposition \ref{Proposition 4.5}, $G_A$ can not be unipotent. It is easy to see that $A$ is not isomorphic to $\mathbb{K}[X]/\langle X\rangle^2 $. Therefore, it is clear that $G_A$ is not nilpotent by Theorem \ref{Theorem 3.6}. Here, $\text{dim}(R/R^2)=3$ and there are no restrictions on the number of generators of $I$ as $l\geq 5$. As a result, we can not use Theorem \ref{Theorem 3.8} and Theorem \ref{Theorem 3.7}. The minimal degree subspace $ W$ associated to $A$ is given by $W=\text{Span}\{X_1^2+X_2^2, X_2^2+X_3^2, X_2X_3\}$. Let $f\in W$, then $f=c_1(X_1^2+X_2^2)+c_2(X_2^2+X_3^2)+c_3(X_2X_3)$ with $c_i\in \mathbb{K}$. Hence, in this example, $W$ does not contain any non-singular homogeneous polynomial of degree at least $3$. However, $\text{dim }W=3=\text{dim}(R/R^2)$ and the ordered sequence $\{X_1^2+X_2^2, X_2^2+X_3^2, X_2X_3\}$ is a regular sequence in $\mathbb{K}[X_1, X_2,X_3]$. Now using the condition $(1)$ of Theorem \ref{Theorem 5.1}, we can conclude that $G_A$ is solvable.
\end{example}
\vskip2mm
\noindent
Similarly, for every $n\geq 3$, we now construct an algebra that satisfies the condition $(1)$ of Theorem \ref{Theorem 5.1}.
\begin{example}\label{Example 5.6}
Let $A_n=\mathbb{K}[X_1,\dots,X_n]/I_n$, where \[I_n=\langle X_1,\dots,X_n\rangle^l+ \langle  X_1^2+X_2^2, X_2^{2}+X_3^{2},\dots, X_{n-1}^{2}+X_{n}^{2}, X_{n-1}X_{n}, X_1^3+X_2^3, X_1^4+X_2^4 \rangle\] with $l\geq 5$, $\text{char }\mathbb{K}=0$ and $n\geq 3$. It is clear that $A_n$ is a non-monomial local graded by the radical algebra for each $n\geq 3$. Similar to the above example, $G_{A_n}$ is not a nilpotent group by Theorem \ref{Theorem 3.6}. Here, $\text{dim}(R_n/R_n^2)=n>2$, and there are no restrictions on the number of generators of $I_n$ for each $n$ as $l\geq 5$. As a result, we can not apply the well-known results (Theorem \ref{Theorem 3.6}, Theorem \ref{Theorem 3.7}, Theorem \ref{Theorem 3.8}, Theorem \ref{Theorem 3.9}). If we calculate the minimal degree subspace $W_n$ associated to $A_n$, then $W_n=\text{ Span}\{X_1^2+X_2^2, X_2^{2}+X_3^{2},\dots, X_{n-1}^{2}+X_{n}^{2}, X_{n-1}X_{n}\}$. It is easy to see that $W_n$ does not contain any non-singular homogeneous polynomial of degree at least $3$ for each $n\geq 3$. However, $\text{dim }W_n=\text{dim}(R_n/R_{n}^2)=n$ and the ordered sequence $\{X_1^2+X_2^2, X_2^{2}+X_3^{2}, \dots,  X_{n-1}^{2}+X_{n}^{2}, X_{n-1}X_{n}\}$ is a regular sequence in $\mathbb{K}[X_1,\dots,X_n]$ for each $n\geq 3$. Therefore, by the condition $(1)$ of Theorem \ref{Theorem 5.1}, $G_{A_n}$ is solvable. The solvability of $G_A$ is a necessary condition for a local finite-dimensional commutative algebra $A$ to occur as the moduli algebra of some isolated hypersurface singularity, as discussed in Section \ref{sec: intro}. Hence, this new class of algebras may occur as the moduli algebra of some isolated hypersurface singularity.
\end{example}
\vskip2mm

\begin{example}\label{Example 5.7}
Let $A=\mathbb{K}[X_1, X_2, X_3]/I$, where \[I=\langle X_1,X_2,X_3\rangle^l+\langle X_1^3+X_1X_2^2+X_2^3,X_2^3+X_3^3, X_1^4+X_2^2X_3^2\rangle\] with $l\geq 5$ and $\text{char }\mathbb{K}=0$. It is clear that $A$ is not a monomial algebra, as not all the generators of $I$ are monomials. Moreover, $A$ is a local graded by the radical algebra as $A/R=\mathbb{K}$ and $I$ is a homogeneous ideal. Since $\mathbb{G}_m$ embeds in $G_A$ by Proposition \ref{Proposition 4.5}, it is not a unipotent group. If $G_A$ is nilpotent, then it must be a torus by Theorem \ref{Theorem 3.6}, which is clearly a contradiction as $A$ is not isomorphic to $\mathbb{K}[X]/\langle X\rangle^2$. Therefore, we get that $G_A$ is not nilpotent. Here, the $\text{dim}(R/R^2)=3=\text{dim }(\langle X_1,X_2,X_3\rangle/\langle X_1,X_2,X_3\rangle^2)>2$. It is also easy to see that there are no restrictions on the number of generators of  $I$ as $l\geq 5$. Therefore, the well-known results do not cover this type of algebra. However, if we calculate the minimal degree subspace $W$ for $A$, then $W=\text{Span}\{X_1^3+X_1X_2^2+X_2^3,X_2^3+X_3^3\}$. It is clear that $f=X_1^3+2X_2^3+X_1X_2^2+X_3^3\in W$, which is a non-singular homogeneous polynomial of degree $3$. Therefore, by Theorem \ref{Theorem 5.1}, $G_A$ is solvable as the condition $(2)$ in Theorem \ref{Theorem 5.1} is satisfied by the minimal degree subspace $W$ associated to $A$.
\end{example}
\vskip2mm
\noindent
Similarly, we now construct an algebra for each $n\geq 3$ which satisfies the condition $(2)$ of Theorem \ref{Theorem 5.1}.
\begin{example}\label{Example 5.8}
Let $A_n=\mathbb{K}[X_1,\dots,X_n]/I_n$, where \[I_n=\langle X_1,\dots,X_n\rangle^l+ \langle  X_1^t+X_2^t, X_2^{t}+X_3^{t}, X_3^{t}+X_4^{t},\dots, X_{n-1}^{t}+X_{n}^{t}, X_1^{t+1}+X_2^{t+1}\rangle\] with $\text{ char }\mathbb{K}=0$, $l\geq (t+2)$, $t\geq 3$ and $n\geq 3$. It is clear that $A_n$ is a local graded by the radical algebra which is not a monomial algebra for each $n\geq 3$. Here, $\text{ dim}(R_n/R_n^2)=n>2$. Since $\mathbb{G}_m$ embeds in $G_{A_n}$ by Proposition \ref{Proposition 4.5}, it is not a unipotent group for each $n\geq 3$. Even $A_n$ is not isomorphic to $\mathbb{K}[X]/\langle X\rangle^2$ for each $n\geq 3$, which implies $G_{A_n}$ is not nilpotent by Theorem \ref{Theorem 3.6}. There are also no restrictions on the number of generators of $I_n$ as $l\geq (t+2)$ for each $n\geq 3$. Therefore, we can not apply any known results (Theorem \ref{Theorem 3.6}, Theorem \ref{Theorem 3.7}, Theorem \ref{Theorem 3.8}, Theorem \ref{Theorem 3.9}). However, if we calculate the minimal degree subspace $W_n$ for $A_n$, then $W_n=\text{Span}\{X_1^t+X_2^{t}, X_2^{t}+X_3^{t}, X_3^{t}+X_4^{t},\dots, X_{n-1}^{t}+X_{n}^{t}\}$. It is easy to see that  $f_n=X_1^{t}+2X_2^{t}+2X_3^{t}+\dots+2X_{n-1}^{t}+X_{n}^{t}\in W_n$ for each $n\geq 3$, which is a non-singular homogeneous polynomial of degree $t\geq 3$. Therefore, by condition $(2)$ of Theorem \ref{Theorem 5.1}, $G_{A_n}$ is solvable. For each $n\geq 3$ and $t\geq 3$, we get an algebra that satisfies the condition $(2)$ of Theorem \ref{Theorem 5.1}.  Hence, it indicates that this new class of algebras captured by the condition $(2)$ in Theorem \ref{Theorem 5.1} may occur as the moduli algebra of some isolated hypersurface singularity, as discussed in the Introduction (Section \ref{sec: intro}).
\end{example}
\vskip2mm
\begin{remark}\label{Remark 5.9}
From Example \ref{Example 5.5}, Example \ref{Example 5.6}, Example \ref{Example 5.7}, and Example \ref{Example 5.8}, it is evident that the new classes of algebras exist that fulfil the conditions of Theorem \ref{Theorem 5.1} and which are not included in well-known results. In addition, by applying the technique described above, it is possible to construct more new local graded by the radical algebras by selecting distinct classes of regular sequences of homogeneous polynomials of degree at least $2$ and non-singular homogeneous polynomials of degree at least $3$, ensuring that $G_A$ is solvable. Motivated by the study of rational Serre spectral sequence, we know that Halperin's conjecture is one of the most prominent open problems in rational homotopy theory. The equivalent algebraic version of this conjecture talks about the non-existence of non-zero negative weight derivations for Artinian complete intersection algebras (which are of the form $\mathbb{C}[X_1,\dots, X_n]/\langle f_1,\dots,f_n\rangle$, where $\{f_1,\dots,f_n\}$ is a regular sequence of weighted homogeneous polynomials in $\mathbb{C}[X_1,\dots, X_n])$ for more detail, see Section $1$ in \cite{Perepechko}, and Section $2$ in \cite{Chen}. The above classes of local algebras are of interest from the perspective of isolated hypersurface singularities and Halperin's conjecture.
\end{remark}
\vskip2mm
\section{Some Remarks and Examples}\label{sec: remark}
\noindent
In this section, we will discuss some examples. We will show that the conditions in Theorem \ref{Theorem 5.1} are sufficient but not necessary.

\vskip2mm
\begin{example}\label{Example 6.1}
  Let $A=\mathbb{K}[X_1,\dots,X_n]/(\langle X_1,\dots, X_n\rangle^l+\langle X_1,\dots ,X_{n-1}\rangle^m)$ with $m<l$ and $n>2$, then $\text{GL}_{n-1}(\mathbb{K})\subset \text{Im}(\Phi_A)$, which implies that $G_A$ is not solvable. It is easy to observe that here $W$ does not contain any non-singular homogeneous polynomial.    
\end{example}
\vskip2mm
\begin{example}\label{Example 6.2}
The conditions in Theorem \ref{Theorem 5.1} are sufficient, but not necessary. Consider the algebra $\mathbb{K}[X,Y]/I$ with $\text{char }\mathbb{K}=0$, where $I=\langle X, Y \rangle^{12}+\langle X^{4}Y^{5}\rangle$. Then, by Theorem \ref{Theorem 3.8}, the group $G_A$ is solvable. However, $W$ does not contain any non-singular homogeneous polynomial. Similarly, if $A=\mathbb{K}[X, Y]/I$ with $\text{char }\mathbb{K}=0$, where $I=\langle X, Y\rangle^5+\langle X^{2}Y, XY^{2}\rangle$, then $W=\text{Span }\{X^{2}Y, XY^{2}\}$ does not have any basis elements which form a regular sequence, but still $G_A$ is solvable by Theorem \ref{Theorem 3.8}.
\end{example}
\vskip2mm
\begin{example}\label{Example 6.3}
Let $A\cong \mathbb{K}[X_1,\dots,X_n]/I$ with $\text{char }\mathbb{K}\neq2$, where $I=\langle X_1,\dots,X_n\rangle^l+\langle q\rangle$, $l>2$, $n>2$ and $q$ is a non-singular quadratic form. Then $G_A\cong \mathbb{G}_m.\text{SO}(n,q)\ltimes U$, by Proposition \ref{Proposition 4.10}. Therefore, $G_A$ is not solvable as $\text{SO}(n,q)$ is not solvable for $n>2$. Here, $W=\text{Span }\{q\}$. This illustration demonstrates that in Theorem \ref{Theorem 5.1}, the condition $(2)$ cannot be relaxed for a non-singular homogeneous polynomial of degree $2$ in $W$.
\end{example}
\vskip2mm
\begin{example}\label{Example 6.4}
    Theorem \ref{Theorem 5.1} does not grant the solvability of $\text{Aut}_{\mathbb{K}}(A)$. For example, take $A=\mathbb{K}[X_1,\dots,X_5]/I$ with $\text{char }\mathbb{K}=0$ or $\text{char }\mathbb{K}>5$, where $I=\langle X_1,\dots,X_5\rangle^l+\langle X_1X_2X_3X_4X_5,X_1^5+X_2^5+\dots X_5^5\rangle$ and $l>6$. Here, \[W=\text{Span }\{X_1X_2X_3X_4X_5, X_1^5+X_2^5+\dots +X_5^5\}.\] By Theorem \ref{Theorem 5.1}, $G_A$ is solvable as $W$ contains a non-singular element, $X_1^5+X_2^5+\dots+X_5^5$ with characteristic of $\mathbb{K}$ is either $0$ or greater than $5$. However, $\text{Aut}_{\mathbb{K}}(A)$ is not solvable as $S_5$ is contained in $\text{Aut}_{\mathbb{K}}(A)$. 
\end{example}
\vskip2mm
\begin{remark}\label{Remark 6.5}
The hypothesis $\text{dim }W=\text{dim}(R/R^2)$ is essential in Theorem \ref{Theorem 5.1}. Let $(G, V)$ be as in Lemma \ref{Lemma 4.12}. Then we can find a homogeneous polynomial $f\in \mathbb{K}[V]^{G}$ such that $\text{Stab}^0(f)=G$, where $G$ is a simple group. Now, consider $A=\mathbb{K}[X_1,\dots,X_n]/I$ with $I=\langle X_1,\dots,X_n\rangle^l+\langle f\rangle$, where $n=\text{dim }V$ and $l>\text{deg}(f)$. Then $G_A\cong \text{Sim}^{0}(f)\ltimes U$ (Proposition \ref{Proposition 4.10}). Therefore, in this situation $\text{dim }W=1<\text{dim}(R/R^2)$. It is easy to see that $f$ is a regular sequence in $\mathbb{K}[X_1,\dots,X_n]$, but $G_A$ is not solvable as it contains a non-trivial simple group $\text{Stab}^{0}(f)$.
\end{remark}
\vskip2mm
\begin{definition}\label{Definition 6.6}
A homogeneous polynomial that is not a monomial is said to be an \emph{$s$-homogeneous polynomial} if it is generated by indeterminants $X_s,\dots, X_n$, $s\geq 1$.
\end{definition}
\vskip2mm
\begin{definition}\label{Definition 6.7}
Let $A$ be a local commutative algebra with Lowey length $l$ and $\text{dim}(R/R^2)=n$ over $\mathbb{K}$. It is proved in Lemma \ref{Lemma 4.7} that the corresponding admissible ideal $I$ for $A$ is generated by all monomials of degree $l$ and non-zero polynomials $P_i$ of degree at least $2$ and at most $l-1$ with no linear component, i.e., $I=\langle X_1,\dots,X_n\rangle^l+\langle P_1,\dots, P_m\rangle$, where $2\leq \text{deg}(P_i)\leq l-1$ and $P_i$ has no linear component for all $i$. The admissible ideal $I$ is said to satisfy \emph{Property $\#$ for $r$}, where $r\in\mathbb{N}$, if $I$ is generated by all monomials of degree $l$ and polynomials $P_i$, where $P_i$ is an $s_i$-homogeneous polynomial for $s_i\geq r$ of degree at least $2$ and at most $l-1$ for all $1\leq i\leq m$.
\end{definition}

\vskip2mm
\begin{remark}\label{Remark 6.8}
Suppose $G_A$ is a solvable group for a graded by the radical local commutative algebra $A$ of $\text{dim }(R/R^2)\geq 3$. In that case, the admissible ideal $I$ corresponding to $A$ with Lowey length $l$ cannot satisfy the Property $\#$ for $r\geq 3$. That is for every generator $f\in I$ with $2\leq \text{deg}(f)\leq l-1$ and no linear component, we can not have $f$ is an $s$-homogeneous polynomial for $s\geq r$. Suppose $I$ satisfies the Property $\#$ for $r\geq 3$, then we have an embedding of $\text{GL}_2(\mathbb{K})$ in $G_A$.  Hence, $G_A$ can not be solvable.
\vskip2mm
\noindent
Suppose the admissible ideal $I$ satisfies the property $\#$ for $r=3$, then the change of variables which is given explicitly by 
\[\biggl\{\Gamma(X_1)=aX_1+bX_2, \Gamma(X_2)=cX_1+dX_2, \Gamma(X_i)=X_i,i\neq 1,2: ad-bc\neq 0\biggl\},\] i.e., \[\Gamma=\begin{pmatrix}
    a & c & 0 \dots 0\\
    b & d & 0 \dots 0\\
    0 & 0 & 1 \dots 0\\
    \vdots & \vdots & 0 \cdots \vdots\\
    0 & 0 & 0 \dots 1
\end{pmatrix}\] 
\vskip2mm
\noindent
stabilizes $I$, hence embeds in $G_A$. Therefore, $\text{GL}_2(\mathbb{K})$ embeds in $G_A$.
\end{remark}
\vskip2mm
\noindent
 We have only given a sufficient condition on $A$ so that $G_A$ is solvable for a local graded by the radical commutative algebra $A$. Still, the following question remains open, originally due to Pollack.
 \begin{question}\label{Question 6.9}
 Let $A$ be a finite-dimensional commutative associative local algebra over an algebraically closed field $\mathbb{K}$. What is a necessary and sufficient condition on the commutative associative algebra $A$, so that $G_A$ is solvable?
 \end{question}
\vskip2mm
\noindent
\textbf{Acknowledgements:}
I would like to thank my supervisor, Prof. Maneesh Thakur, for valuable discussions with me during the preparation of this article and for giving some fruitful suggestions on the preliminary version. I want to thank my friend Sayan Pal for making some useful suggestions. I am also grateful to the Indian Statistical Institute, Bangalore for providing
an excellent and productive research atmosphere. Finally, I thank the referee for the very
useful comments and suggestions, which have improved the readability of the paper.

\bibliography{references}

@article {Pollack,
    AUTHOR = {Pollack, R. David},
     TITLE = {Algebras and their automorphism groups},
   JOURNAL = {Comm. Algebra},
  FJOURNAL = {Communications in Algebra},
    VOLUME = {17},
      YEAR = {1989},
    NUMBER = {8},
     PAGES = {1843--1866},
      ISSN = {0092-7872,1532-4125},
   MRCLASS = {16A72 (16A46)},
  MRNUMBER = {1013471},
MRREVIEWER = {Susan\ Montgomery},
       DOI = {10.1080/00927878908823824},
       URL = {https://doi-org.libraryisikolkata.remotexs.in/10.1080/00927878908823824},
}

@article {Yau,
    AUTHOR = {Yau, Stephen S.-T.},
     TITLE = {Solvable {L}ie algebras and generalized {C}artan matrices
              arising from isolated singularities},
   JOURNAL = {Math. Z.},
  FJOURNAL = {Mathematische Zeitschrift},
    VOLUME = {191},
      YEAR = {1986},
    NUMBER = {4},
     PAGES = {489--506},
      ISSN = {0025-5874,1432-1823},
   MRCLASS = {32B30 (14B05 17B67)},
  MRNUMBER = {832806},
MRREVIEWER = {Jonathan\ M.\ Wahl},
       DOI = {10.1007/BF01162338},
       URL = {https://doi-org.libraryisikolkata.remotexs.in/10.1007/BF01162338},
}

@article {AS1,
    AUTHOR = {Guil Asensio, Francisco and Saor\'in, Manuel},
     TITLE = {The group of automorphisms of a commutative algebra},
   JOURNAL = {Math. Z.},
  FJOURNAL = {Mathematische Zeitschrift},
    VOLUME = {219},
      YEAR = {1995},
    NUMBER = {1},
     PAGES = {31--48},
      ISSN = {0025-5874,1432-1823},
   MRCLASS = {13E10 (13B10 20G15)},
  MRNUMBER = {1340847},
MRREVIEWER = {James\ K.\ Deveney},
       DOI = {10.1007/BF02572348},
       URL = {https://doi-org.libraryisikolkata.remotexs.in/10.1007/BF02572348},
}

@article {AS2,
    AUTHOR = {Guil-Asensio, Francisco and Saor\'in, Manuel},
     TITLE = {The automorphism group and the {P}icard group of a monomial
              algebra},
   JOURNAL = {Comm. Algebra},
  FJOURNAL = {Communications in Algebra},
    VOLUME = {27},
      YEAR = {1999},
    NUMBER = {2},
     PAGES = {857--887},
      ISSN = {0092-7872,1532-4125},
   MRCLASS = {16W30},
  MRNUMBER = {1672003},
MRREVIEWER = {Stefaan\ Caenepeel},
       DOI = {10.1080/00927879908826466},
       URL = {https://doi-org.libraryisikolkata.remotexs.in/10.1080/00927879908826466},
}

@article {AS3,
    AUTHOR = {Guil-Asensio, Francisco and Saor\'in, Manuel},
     TITLE = {The group of outer automorphisms and the {P}icard group of an
              algebra},
   JOURNAL = {Algebr. Represent. Theory},
  FJOURNAL = {Algebras and Representation Theory},
    VOLUME = {2},
      YEAR = {1999},
    NUMBER = {4},
     PAGES = {313--330},
      ISSN = {1386-923X,1572-9079},
   MRCLASS = {16W20 (16W22)},
  MRNUMBER = {1733381},
MRREVIEWER = {Wenxue\ Huang},
       DOI = {10.1023/A:1009973319703},
       URL = {https://doi-org.libraryisikolkata.remotexs.in/10.1023/A:1009973319703},
}

@article {Perepechko,
    AUTHOR = {Perepechko, Alexander},
     TITLE = {On solvability of the automorphism group of a
              finite-dimensional algebra},
   JOURNAL = {J. Algebra},
  FJOURNAL = {Journal of Algebra},
    VOLUME = {403},
      YEAR = {2014},
     PAGES = {445--458},
      ISSN = {0021-8693,1090-266X},
   MRCLASS = {13N15 (17B40)},
  MRNUMBER = {3166084},
MRREVIEWER = {Michael\ Gr.\ Voskoglou},
       DOI = {10.1016/j.jalgebra.2014.01.018},
       URL = {https://doi-org.libraryisikolkata.remotexs.in/10.1016/j.jalgebra.2014.01.018},
}

@article {GG,
    AUTHOR = {Garibaldi, Skip and Guralnick, Robert M.},
     TITLE = {Simple groups stabilizing polynomials},
   JOURNAL = {Forum Math. Pi},
  FJOURNAL = {Forum of Mathematics. Pi},
    VOLUME = {3},
      YEAR = {2015},
     PAGES = {e3, 41},
      ISSN = {2050-5086},
   MRCLASS = {20G15 (15A72 20G41)},
  MRNUMBER = {3406824},
MRREVIEWER = {Timothy\ C.\ Burness},
       DOI = {10.1017/fmp.2015.3},
       URL = {https://doi-org.libraryisikolkata.remotexs.in/10.1017/fmp.2015.3},
}

@article {Sei,
    AUTHOR = {Seitz, Gary M.},
     TITLE = {The maximal subgroups of classical algebraic groups},
   JOURNAL = {Mem. Amer. Math. Soc.},
  FJOURNAL = {Memoirs of the American Mathematical Society},
    VOLUME = {67},
      YEAR = {1987},
    NUMBER = {365},
     PAGES = {iv+286},
      ISSN = {0065-9266,1947-6221},
   MRCLASS = {20G15 (20E28 20G05)},
  MRNUMBER = {888704},
MRREVIEWER = {James\ E.\ Humphreys},
       DOI = {10.1090/memo/0365},
       URL = {https://doi-org.libraryisikolkata.remotexs.in/10.1090/memo/0365},
}

@book {Curtis,
    AUTHOR = {Curtis, Charles W. and Reiner, Irving},
     TITLE = {Representation theory of finite groups and associative
              algebras},
      NOTE = {Reprint of the 1962 original},
 PUBLISHER = {AMS Chelsea Publishing, Providence, RI},
      YEAR = {2006},
     PAGES = {xiv+689},
      ISBN = {0-8218-4066-5},
   MRCLASS = {16-02 (20-02 20Cxx)},
  MRNUMBER = {2215618},
       DOI = {10.1090/chel/356},
       URL = {https://doi-org.libraryisikolkata.remotexs.in/10.1090/chel/356},
}

@book {Humphreys,
    AUTHOR = {Humphreys, James E.},
     TITLE = {Linear algebraic groups},
    SERIES = {Graduate Texts in Mathematics},
    VOLUME = {No. 21},
 PUBLISHER = {Springer-Verlag, New York-Heidelberg},
      YEAR = {1975},
     PAGES = {xiv+247},
   MRCLASS = {20GXX (14LXX)},
  MRNUMBER = {396773},
MRREVIEWER = {T.\ Ono},
}

@book {Borel,
    AUTHOR = {Borel, Armand},
     TITLE = {Linear algebraic groups},
    SERIES = {Graduate Texts in Mathematics},
    VOLUME = {126},
   EDITION = {Second},
 PUBLISHER = {Springer-Verlag, New York},
      YEAR = {1991},
     PAGES = {xii+288},
      ISBN = {0-387-97370-2},
   MRCLASS = {20-01 (20Gxx)},
  MRNUMBER = {1102012},
MRREVIEWER = {F.\ D.\ Veldkamp},
       DOI = {10.1007/978-1-4612-0941-6},
       URL = {https://doi-org.libraryisikolkata.remotexs.in/10.1007/978-1-4612-0941-6},
}

@book {Pierce,
    AUTHOR = {Pierce, Richard S.},
     TITLE = {Associative algebras},
    SERIES = {Studies in the History of Modern Science},
    VOLUME = {9},
      NOTE = {Graduate Texts in Mathematics, 88},
 PUBLISHER = {Springer-Verlag, New York-Berlin},
      YEAR = {1982},
     PAGES = {xii+436},
      ISBN = {0-387-90693-2},
   MRCLASS = {16-01 (12-01)},
  MRNUMBER = {674652},
MRREVIEWER = {S.\ S.\ Page},
}

@book {Auslander,
    AUTHOR = {Auslander, Maurice and Reiten, Idun and Smal\o, Sverre O.},
     TITLE = {Representation theory of {A}rtin algebras},
    SERIES = {Cambridge Studies in Advanced Mathematics},
    VOLUME = {36},
      NOTE = {Corrected reprint of the 1995 original},
 PUBLISHER = {Cambridge University Press, Cambridge},
      YEAR = {1997},
     PAGES = {xiv+425},
      ISBN = {0-521-41134-3; 0-521-59923-7},
   MRCLASS = {16Gxx (16-02)},
  MRNUMBER = {1476671},
}

@conference{Saorin,
    author ={M. Saorin} ,
    booktitle ={In: J, Abrams, et al: Methods in Module Theory},
    title = {Gradability of algebras},
    year = {1991}
}

@article {Schneider,
    AUTHOR = {Schneider, J. E.},
     TITLE = {Orthogonal groups of nonsingular forms of higher degree},
   JOURNAL = {J. Algebra},
  FJOURNAL = {Journal of Algebra},
    VOLUME = {27},
      YEAR = {1973},
     PAGES = {112--116},
      ISSN = {0021-8693},
   MRCLASS = {20G15},
  MRNUMBER = {325795},
MRREVIEWER = {A.\ O.\ Morris},
       DOI = {10.1016/0021-8693(73)90166-X},
       URL = {https://doi-org.libraryisikolkata.remotexs.in/10.1016/0021-8693(73)90166-X},
}

@book {GKZ,
    AUTHOR = {Gelfand, I. M. and Kapranov, M. M. and Zelevinsky, A.
              V.},
     TITLE = {Discriminants, resultants, and multidimensional determinants},
    SERIES = {Mathematics: Theory \& Applications},
 PUBLISHER = {Birkh\"auser Boston, Inc., Boston, MA},
      YEAR = {1994},
     PAGES = {x+523},
      ISBN = {0-8176-3660-9},
   MRCLASS = {14N05 (13D25 14M25 15A69 33C70 52B20)},
  MRNUMBER = {1264417},
MRREVIEWER = {I.\ Dolgachev},
       DOI = {10.1007/978-0-8176-4771-1},
       URL = {https://doi-org.libraryisikolkata.remotexs.in/10.1007/978-0-8176-4771-1},
}

@article {OS,
    AUTHOR = {Orlik, Peter and Solomon, Louis},
     TITLE = {Singularities. {I}. {H}ypersurfaces with an isolated
              singularity},
   JOURNAL = {Advances in Math.},
  FJOURNAL = {Advances in Mathematics},
    VOLUME = {27},
      YEAR = {1978},
    NUMBER = {3},
     PAGES = {256--272},
      ISSN = {0001-8708},
   MRCLASS = {14B05 (14B15 32C40)},
  MRNUMBER = {476734},
MRREVIEWER = {Peter\ Giblin},
       DOI = {10.1016/0001-8708(78)90101-9},
       URL = {https://doi-org.libraryisikolkata.remotexs.in/10.1016/0001-8708(78)90101-9},
}

@article {Kraft,
    AUTHOR = {Kraft, H. and Procesi, C.},
     TITLE = {Graded morphisms of {$G$}-modules},
   JOURNAL = {Ann. Inst. Fourier (Grenoble)},
  FJOURNAL = {Universit\'e{} de Grenoble. Annales de l'Institut Fourier},
    VOLUME = {37},
      YEAR = {1987},
    NUMBER = {4},
     PAGES = {161--166},
      ISSN = {0373-0956,1777-5310},
   MRCLASS = {20G05 (13F20)},
  MRNUMBER = {927395},
MRREVIEWER = {P.\ E.\ Newstead},
       DOI = {10.5802/aif.1115},
       URL = {https://doi-org.libraryisikolkata.remotexs.in/10.5802/aif.1115},
}

@article {Mather,
    AUTHOR = {Mather, John N. and Yau, Stephen S. T.},
     TITLE = {Classification of isolated hypersurface singularities by their
              moduli algebras},
   JOURNAL = {Invent. Math.},
  FJOURNAL = {Inventiones Mathematicae},
    VOLUME = {69},
      YEAR = {1982},
    NUMBER = {2},
     PAGES = {243--251},
      ISSN = {0020-9910,1432-1297},
   MRCLASS = {32B30 (14B05)},
  MRNUMBER = {674404},
MRREVIEWER = {Alexandru\ Dimca},
       DOI = {10.1007/BF01399504},
       URL = {https://doi-org.libraryisikolkata.remotexs.in/10.1007/BF01399504},
}

@book {Greuel,
    AUTHOR = {Greuel, G.-M. and Lossen, C. and Shustin, E.},
     TITLE = {Introduction to singularities and deformations},
    SERIES = {Springer Monographs in Mathematics},
 PUBLISHER = {Springer, Berlin},
      YEAR = {2007},
     PAGES = {xii+471},
      ISBN = {978-3-540-28380-5; 3-540-28380-3},
   MRCLASS = {32Sxx (14B05)},
  MRNUMBER = {2290112},
MRREVIEWER = {Vasile\ Br\^inz\u anescu},
}

@article {Popov,
    AUTHOR = {Popov, V. L.},
     TITLE = {Group varieties and group structures},
   JOURNAL = {Izv. Ross. Akad. Nauk Ser. Mat.},
  FJOURNAL = {Izvestiya Rossiiskoi Akademii Nauk. Seriya Matematicheskaya},
    VOLUME = {86},
      YEAR = {2022},
    NUMBER = {5},
     PAGES = {73--96},
      ISSN = {1607-0046,2587-5906},
   MRCLASS = {14L15 (14M99 20G20)},
  MRNUMBER = {4582538},
MRREVIEWER = {Nicolae\ Manolache},
       DOI = {10.4213/im9272},
       URL = {https://doi-org.libraryisikolkata.remotexs.in/10.4213/im9272},
}

@article {SK,
    AUTHOR = {Sato, M. and Kimura, T.},
     TITLE = {A classification of irreducible prehomogeneous vector spaces
              and their relative invariants},
   JOURNAL = {Nagoya Math. J.},
  FJOURNAL = {Nagoya Mathematical Journal},
    VOLUME = {65},
      YEAR = {1977},
     PAGES = {1--155},
      ISSN = {0027-7630,2152-6842},
   MRCLASS = {32M10 (20G05)},
  MRNUMBER = {430336},
MRREVIEWER = {A.\ L.\ Onishchik},
       URL = {http://projecteuclid.org/euclid.nmj/1118796150},
}

@article {FR,
    AUTHOR = {Eisele, Florian and Raedschelders, Theo},
     TITLE = {On solvability of the first {H}ochschild cohomology of a
              finite-dimensional algebra},
   JOURNAL = {Trans. Amer. Math. Soc.},
  FJOURNAL = {Transactions of the American Mathematical Society},
    VOLUME = {373},
      YEAR = {2020},
    NUMBER = {11},
     PAGES = {7607--7638},
      ISSN = {0002-9947,1088-6850},
   MRCLASS = {16E40 (16G10 16G60)},
  MRNUMBER = {4169669},
MRREVIEWER = {Leandro\ R.\ Cagliero},
       DOI = {10.1090/tran/8064},
       URL = {https://doi-org.libraryisikolkata.remotexs.in/10.1090/tran/8064},
}

@article {CSS,
    AUTHOR = {Chaparro, Cristian and Schroll, Sibylle and Solotar, Andrea},
     TITLE = {On the {L}ie algebra structure of the first {H}ochschild
              cohomology of gentle algebras and {B}rauer graph algebras},
   JOURNAL = {J. Algebra},
  FJOURNAL = {Journal of Algebra},
    VOLUME = {558},
      YEAR = {2020},
     PAGES = {293--326},
      ISSN = {0021-8693,1090-266X},
   MRCLASS = {16E40 (16G20 16W25)},
  MRNUMBER = {4102135},
MRREVIEWER = {Karin\ Erdmann},
       DOI = {10.1016/j.jalgebra.2020.02.003},
       URL = {https://doi-org.libraryisikolkata.remotexs.in/10.1016/j.jalgebra.2020.02.003},
}

@article {DSSS,
    AUTHOR = {Rubio y Degrassi, Lleonard and Schroll, Sibylle and Solotar,
              Andrea},
     TITLE = {The first {H}ochschild cohomology as a {L}ie algebra},
   JOURNAL = {Quaest. Math.},
  FJOURNAL = {Quaestiones Mathematicae. Journal of the South African
              Mathematical Society},
    VOLUME = {46},
      YEAR = {2023},
    NUMBER = {9},
     PAGES = {1955--1980},
      ISSN = {1607-3606,1727-933X},
   MRCLASS = {16E40 (16D90 16G60)},
  MRNUMBER = {4634188},
MRREVIEWER = {Shengyong\ Pan},
       DOI = {10.2989/16073606.2022.2115424},
       URL = {https://doi-org.libraryisikolkata.remotexs.in/10.2989/16073606.2022.2115424},
}

@article {S,
    AUTHOR = {Strametz, Claudia},
     TITLE = {The {L}ie algebra structure on the first {H}ochschild
              cohomology group of a monomial algebra},
   JOURNAL = {J. Algebra Appl.},
  FJOURNAL = {Journal of Algebra and its Applications},
    VOLUME = {5},
      YEAR = {2006},
    NUMBER = {3},
     PAGES = {245--270},
      ISSN = {0219-4988,1793-6829},
   MRCLASS = {16E40 (17B60)},
  MRNUMBER = {2235810},
MRREVIEWER = {Sarah\ J.\ Witherspoon},
       DOI = {10.1142/S0219498806001120},
       URL = {https://doi-org.libraryisikolkata.remotexs.in/10.1142/S0219498806001120},
}

@incollection {PG,
    AUTHOR = {Gabriel, Peter},
     TITLE = {Auslander-{R}eiten sequences and representation-finite
              algebras},
 BOOKTITLE = {Representation theory, {I} ({P}roc. {W}orkshop, {C}arleton
              {U}niv., {O}ttawa, {O}nt., 1979)},
    SERIES = {Lecture Notes in Math.},
    VOLUME = {831},
     PAGES = {1--71},
 PUBLISHER = {Springer, Berlin},
      YEAR = {1980},
      ISBN = {3-540-10263-9},
   MRCLASS = {16A64 (16A46)},
  MRNUMBER = {607140},
MRREVIEWER = {Idun\ Reiten},
}

@article {YAU2,
    AUTHOR = {Yau, Stephen S.-T.},
     TITLE = {Solvability of {L}ie algebras arising from isolated
              singularities and nonisolatedness of singularities defined by
              {${\rm sl}(2,{\bf C})$} invariant polynomials},
   JOURNAL = {Amer. J. Math.},
  FJOURNAL = {American Journal of Mathematics},
    VOLUME = {113},
      YEAR = {1991},
    NUMBER = {5},
     PAGES = {773--778},
      ISSN = {0002-9327,1080-6377},
   MRCLASS = {32S25 (17B30 32S05)},
  MRNUMBER = {1129292},
MRREVIEWER = {Aleksandr\ G.\ Aleksandrov},
       DOI = {10.2307/2374785},
       URL = {https://doi-org.libraryisikolkata.remotexs.in/10.2307/2374785},
}

@article{dd,
      title={{On $R$-equivalence of automorphism groups of associative algebras}}, 
      author={Dibyendu Das},
      year={2025},
      note={\url{https://arxiv.org/abs/2512.24357}}, 
}

@book {ISD,
    AUTHOR = {Assem, Ibrahim and Simson, Daniel and Skowro\'nski, Andrzej},
     TITLE = {Elements of the representation theory of associative algebras.
              {V}ol. 1},
    SERIES = {London Mathematical Society Student Texts},
    VOLUME = {65},
      NOTE = {Techniques of representation theory},
 PUBLISHER = {Cambridge University Press, Cambridge},
      YEAR = {2006},
     PAGES = {x+458},
      ISBN = {978-0-521-58423-4; 978-0-521-58631-3; 0-521-58631-3},
   MRCLASS = {16G10 (16-02)},
  MRNUMBER = {2197389},
MRREVIEWER = {Peter\ W.\ Donovan},
       DOI = {10.1017/CBO9780511614309},
       URL = {https://doi-org.libraryisikolkata.remotexs.in/10.1017/CBO9780511614309},
}

@article {Chen,
    AUTHOR = {Chen, Hao and Yau, Stephen S.-T. and Zuo, Huaiqing},
     TITLE = {Non-existence of negative weight derivations on positively
              graded {A}rtinian algebras},
   JOURNAL = {Trans. Amer. Math. Soc.},
  FJOURNAL = {Transactions of the American Mathematical Society},
    VOLUME = {372},
      YEAR = {2019},
    NUMBER = {4},
     PAGES = {2493--2535},
      ISSN = {0002-9947,1088-6850},
   MRCLASS = {13N15 (14B05 32S05)},
  MRNUMBER = {3988584},
MRREVIEWER = {Kiyoshi\ Baba},
       DOI = {10.1090/tran/7628},
       URL = {https://doi.org/10.1090/tran/7628},
}
\bibliographystyle{amsplain}
\end{document}